\documentclass[11pt]{article}

\usepackage{amssymb,latexsym}
\usepackage{pb-diagram}
\usepackage{amsmath,amscd}
\usepackage{theorem} 
\usepackage[all]{xy}

\title{\bf Picard groups in Poisson geometry}
\author{Henrique Bursztyn\thanks{{\tt henrique@math.toronto.edu}}\\[0.1cm]
        Department of Mathematics\\
	University of Toronto\\
        Toronto, Ontario M5S 3G3, Canada
        \\[0.2cm]
        Alan Weinstein \thanks{{\tt alanw@math.berkeley.edu}, Research
	partially supported by NSF grants DMS-9971505 and DMS-0204100} \\[0.1cm]
         Department of Mathematics\\
        University of California, Berkeley\\
        CA, 94720-3840}
\date{} 

\newcommand{\Chi}        {{\mathcal X}}

\newcommand{\id}         {{\mathrm {Id}}}

\newcommand{\Ker}        {{\mathrm {Ker}}}

\newcommand{\Aut}        {{\mathrm {Aut}}}
\newcommand{\Inaut}      {{\mathrm {InnAut}}}
\newcommand{\Outaut}     {{\mathrm {OutAut}}}

\newcommand{\pr}         {{\mathrm{pr}}}
\newcommand{\Diff}       {\mathrm{Diff}}
\newcommand{\Poiss}      {\mathrm{Aut}}
\newcommand{\InnPois}   {\mathrm{InnAut}}
\newcommand{\OutPois}   {\mathrm{OutAut}}
\newcommand{\Sympl}      {\mathrm{Sympl}}

\newcommand{\st} [1]     {\scriptscriptstyle{{#1}}}
\newcommand{\univ} [1]   {\widetilde{{#1}}}

\newcommand{\inc}        {j}
\newcommand{\cent}       {\mathcal{Z}}
\newcommand{\proj}       {h}
\newcommand{\LBis}       {\mathrm{LBis}}
\newcommand{\Bis}        {\mathrm{Bis}}
\newcommand{\IsoLBis}    {\mathrm{IsoLBis}}
\newcommand{\CIsoBis}   {\mathrm{ZIsoBis}}

\newcommand{\frakt}     {\mathfrak{t}}
\newcommand{\gstar}     {\mathfrak{g}^*}

\newcommand{\catalg}        {\mathsf{Alg}}
\newcommand{\catgr}        {\mathsf{LG}} 
\newcommand{\catsym}        {\mathsf{SG}} 
\newcommand{\catpoiss}        {\mathsf{Poiss}} 

\newcommand{\A}          {\mathcal{A}}
\newcommand{\B}          {\mathcal{B}}
\newcommand{\C}          {\mathcal{C}}
\newcommand{\D}          {\mathcal{D}}
\newcommand{\calo}          {\mathcal{O}}

\newcommand{\Pic}         {\mathrm{Pic}}
\newcommand{\Piclie}      {\mathfrak{pic}}

\newcommand{\PicA}        {\mathrm{Pic}_{\st{\A}}}
\newcommand{\PicZ}        {\mathrm{Pic}_{\st{\cent(\A)}}}
\newcommand{\PicZP}       {\mathrm{Pic}_{\st{\cent(P)}}}
\newcommand{\PicZG}       {\mathrm{Pic}_{\st{\cent(\Gamma)}}}

\newcommand{\X}          {\mathrm{X}}
\newcommand{\Y}          {\mathrm{Y}}
\newcommand{\Z}          {\mathrm{Z}}

\newcommand{\AXB}        {_{\scriptscriptstyle{\A}}{\X}_{\scriptscriptstyle{\B}}}

\newcommand{\BYC}        {_{\scriptscriptstyle{\B}}{\Y}_{\scriptscriptstyle{\C}}}

\newcommand{\CZD}        {_{\scriptscriptstyle{\C}}{\Z}_{\scriptscriptstyle{\D}}}
\newcommand{\PSPn}        {P\stackrel{J_1}{\leftarrow}S\stackrel{J_2}{\rightarrow}P}
\newcommand{\PSPb}        {P_1\stackrel{J_1}{\leftarrow}S\stackrel{J_2}{\rightarrow}\overline{P}_2}

\newcommand{\PSP}        {P_1\stackrel{J_1}{\leftarrow}S\stackrel{J_2}{\rightarrow}P_2}
\newcommand{\PSPprime} {P_2\stackrel{J_2'}{\leftarrow}S'\stackrel{J_3'}{\rightarrow}P_3}

\newcommand{\tensorB}    {\otimes_{\scriptscriptstyle{\B}}}

\newcommand{\reals}      {{\mathbb R}}

\newcommand{\complex}    {{\mathbb C}}

\newtheorem{lemma} {Lemma} [section]
\newtheorem{proposition} [lemma] {Proposition}

\newtheorem{theorem} [lemma] {Theorem}
\newtheorem{corollary} [lemma] {Corollary}

\theorembodyfont{\rm} 

\newtheorem{example}[lemma] {Example}

\newtheorem{remark}[lemma]{Remark}

\newenvironment{proof}{{\sc Proof:}}{{\hspace*{\fill} $\square$\\}}

\numberwithin{equation}{section}

\begin{document}

\maketitle

\begin{abstract}
We study isomorphism classes of symplectic dual pairs $P\leftarrow
S\rightarrow\overline{P}$, where $P$ is an integrable Poisson
manifold, $S$ is symplectic, and the two maps are complete, surjective
Poisson submersions with connected and simply-connected fibres.  
For fixed $P$, these Morita
self-equivalences of $P$ form a group
$\Pic(P)$ under a natural ``tensor product'' operation.  Variants of
this construction are also studied, for rings (the origin of the
notion of Picard group), Lie groupoids, and symplectic groupoids.
\end{abstract}

\begin{center}
{\bf DEDICATED TO PIERRE CARTIER}
\end{center}

\section{Introduction}

In his article on the occasion of the 40th anniversary of the IHES,
Pierre Cartier \cite{Ca2001} stressed the idea of
``representations as points,''  which was central to the work
of Grothendieck,  and which leads to categories in which 
 morphisms map representations rather than
``ordinary'' points.  In algebra, this idea underlies
the notion
of {\bf Morita equivalence} \cite{Morita58}, 
and the same term is now applied
to objects which ``encode'' algebras, such
as groupoids and Poisson manifolds.  In all these
settings, a morphism between objects $A$ and $B$ is an isomorphism class
of ``bimodules'' of some sort, on which $A$ and $B$ have commuting actions from
the left and right respectively.  Morphisms are composed by an
associative ``tensor product'' operation, and the resulting group of
self-equivalences of an object is called its {\bf Picard group}.  

Morita equivalence for Poisson manifolds and 
symplectic groupoids was introduced by Xu \cite{Xu91,Xu91b}, based on
notions of equivalence of topological groupoids by Haefliger
\cite{Ha84}, Muhly, Renault, and Williams \cite{MuReWi87}, Moerdijk 
\cite{Moerd91} and
others, as well as on the idea of a symplectic dual pair \cite{We83}.
The latter was itself an extension of Howe's notion of linear dual
pair  \cite{Ho79}, also introduced 
for the purpose of obtaining equivalences of representation categories.
Xu's ideas and their relation to Rieffel's \cite{Ri74}
strong Morita equivalence of $C^*$-algebras have been explained and
extended by Landsman in a series of publications
\cite{Land98,Land00b,Land01}, of which the present work may be seen
as a further extension.

In this paper, we briefly review the work of Xu and Landsman, and then
study several interesting examples, concentrating on their  Picard
groups.  We also discuss briefly another notion of
equivalence in Poisson geometry, {\bf gauge equivalence}
\cite{SeWe01}, which has been shown \cite{BuWa2001,JSW2001} to be
related to algebraic Morita equivalence via deformation quantization;
we revisit and extend the results in \cite{BuRad02} relating gauge
equivalence to geometric Morita equivalence.  Work of Rieffel and
Schwarz \cite{RiefSch} also suggests a relation between gauge
equivalence and strong Morita equivalence of $C^*$-algebras; also see
\cite{TaWe}, in which Dirac structures play a key role.

In the final sections, we briefly treat some subjects which we have
just begun to study.  We investigate the Lie algebras of Picard groups
in various categories.   We
also suggest a way of ``enriching'' the structure of the category of
representations of a Poisson manifold so that representation
equivalence implies Morita equivalence.  Finally, we
suggest a way of extending
the notion of Morita equivalence for
groupoids and Poisson manifolds.  By enlarging the class of bimodules
from manifolds to more general leaf spaces of foliations, we hope to
eliminate some technical restrictions (integrability and regularity)
which have limited previous work on geometric Morita equivalence.
Here, the bimodules themselves are Morita equivalence classes! 

The paper as a whole is somewhat like \cite{We97} in that we introduce a new
concept in Poisson geometry, study some examples, and ask a lot of
questions.  

\bigskip

\noindent
{\bf Acknowledgments.}  This paper is a substantially expanded
version of a talk entitled ``Notions d'\'equivalence en g\'eom\'etrie de
Poisson,'' which was presented by A.W. at the Colloque en l'honneur de
Pierre Cartier in June 2002.  He would like to thank the
organizers for the invitation, which stimulated some of the work
described here.  
H.B. thanks the Universit\'e Libre de Bruxelles and the Institute for Pure and
Applied Mathematics (UCLA) for their hospitality while part of 
this work was being done.
We would also like to thank Bill Arveson, Marius Crainic,
Hans Duistermaat, Bert Kostant, Klaas Landsman, Peter Michor, Ieke
Moerdijk, Janez Mr\v{c}un, Tudor Ratiu, Dmitry Roytenberg,
Saul Schleimer, and Ping Xu for helpful
advice.  Finally, we give special thanks to Olga Radko for her
substantial contribution to the writing of Section \ref{sec:TSS1}.

\section{Bimodules as generalized morphisms of algebras}

Before passing to geometry, we will review the idea of generalized
morphisms in the algebraic context where it first arose.

Let us consider the category whose objects are  unital algebras 
over a fixed (commutative, unital) ring $k$
and arrows are algebra homomorphisms.  Since homomorphisms act from
the left, we adopt the convention that a homomorphism $\phi$ from $\B$ to
$\A$ is generally written as an arrow $\A \stackrel{\phi}{\leftarrow} \B$.
We will denote the automorphism group of a $k$-algebra $\A$ by $\Aut(\A)$.

We can define another category with the same objects,
 a ``generalized morphism'' $\A \leftarrow \B$ being an $(\A,\B)$-bimodule
$\AXB$, i.e. a $k$-module which is a 
left $\A$-module and right $\B$-module on which the
actions of the two algebras respect the $k$-module structure and commute. 
The ``composition'' of $\AXB$ and $\BYC$ is 
defined to be their tensor product over $\B$, an $(\A,\C)$-bimodule.
Any algebra homomorphism $ A \stackrel{\phi}{\leftarrow} \B$ can be regarded
as such a generalized morphism, namely
the $(\A,\B)$-bimodule $\A_{\st{\phi}}$, with
$\A$-action given by left multiplication
and right $\B$-action defined by $a \cdot b := a \phi(b)$.  

A feature of the more general 
``morphisms'' (and the reason for the quotes) is that, strictly speaking,
they do not define a category.  Rather, they are 
horizontal arrows in  a {\it bicategory}, see e.g. \cite{Ben67,Mclan71}.
For example, their composition is not associative, but just
associative up to a (natural) 
bimodule isomorphism:
$$
({\AXB} \tensorB  {\BYC})\otimes_{\st{\mathcal C}}{\CZD}
\cong
{\AXB} \tensorB ({\BYC}\otimes_{\st{\mathcal C}}{\CZD}).
$$
Similarly, if we have homomorphisms $ \A \stackrel{\phi}{\leftarrow}
\B$ and 
$\B \stackrel{\psi}{\leftarrow} \C$, then the $(\A,\C)$-bimodules $\A_{\phi\psi}$ and
$\A_\phi \tensorB \B_\psi$ are naturally isomorphic, but not
equal.  We therefore introduce a category $\catalg$ in which the
objects are unital algebras (the base ring $k$ being fixed) and the morphisms $\A
\leftarrow \B$ are
isomorphism classes of $(\A,\B)$-bimodules.  The map
$\phi\mapsto \A_{\st{\phi}}$ induces a functor $\inc$ from the usual
category of algebras to $\catalg$.  We will soon see that $\inc$ is
not faithful.

 Invertible morphisms in
$\catalg$ are known as {\bf Morita equivalences}, and the group of
automorphisms of $\A$ in $\catalg$ is called its
{\bf Picard group}, denoted by $\Pic(\A)$. 
The functor $\inc$ restricts to a group homomorphism from $\Aut(\A)$
to $\Pic(\A)$.  

A simple computation shows that the automorphisms which become trivial
in the Picard group are just the inner automorphisms; i.e., we have the
following exact sequence of groups \cite{Bass68}:
\begin{equation}\label{eq:exactsq}
1 \rightarrow \Inaut(\A) \rightarrow \Aut(\A) \stackrel{\inc}{\rightarrow} \Pic(\A).
\end{equation}
In particular, when $\A$ is commutative, $\Aut(\A)$ sits inside $\Pic(\A)$ as a subgroup
via the map $\inc$.

Now let $\cent(\A)$ denote the center of $\A$. There is a
natural group homomorphism
\cite{Bass68}
\begin{equation}\label{eq:hhom}
\proj: \Pic(\A) \rightarrow \Aut(\cent(\A))
\end{equation}
which takes (the isomorphism class of) each invertible $(\A,\A)$-bimodule $\X$ to 
$\proj_{\X} \in \Aut(\cent(\A))$, defined by the condition that
$\proj_{\X}(z)x = xz$ for all $z \in \cent(\A)$ and $x \in \X$.

If we denote by $\PicZ(\A)$ the subgroup of $\Pic(\A)$ given by
bimodules $\X$
satisfying $zx =xz$ for all $x \in \X$ and $z \in \cent(\A)$, then we have
the following exact sequence:
\begin{equation}\label{eq:exactsq2}
1 \rightarrow \PicZ(\A) \rightarrow \Pic(\A) \stackrel{\proj}{\rightarrow} \Aut(\cent(\A)).
\end{equation}

When $\A$ is commutative, a simple computation shows that
the map $\proj$ is split by the map $\inc$.
In this case, one can write $\Pic(\A)$ as a semi-direct product of $\Aut(\A)$
and the group
$\PicA(\A):= \PicZ(\A)$:
$$
\Pic(\A) = \Aut(\A) \ltimes \PicA(\A).
$$ 
The action of $\Aut(\A)$ on $\PicA(\A)$ is given by
$\X \stackrel{\phi}{\mapsto} {_{\st{\phi}}}\X_{\st{\phi}}$, where
the left and right $\A$-module structures on ${_{\st{\phi}}}\X_{\st{\phi}}$
are given by $a \cdot x := \phi(a)x$ and $x \cdot a := x \phi(a)$.
The group $\PicA(\A)$, consisting of bimodules for which the left and
right module structures are the same, is often called the {\it commutative} Picard
group of $\A$ (to avoid confusion with $\Pic(\A)$, which makes sense
even when $\A$ is noncommutative).

\begin{example}\label{ex:linebundles}
Let $M$ be a smooth manifold, and let $\A =  C^\infty(M)$ be the
${\mathbb C}$-algebra of complex-valued 
smooth functions on $M$. 
In this case,
\begin{equation}\label{eq:identif}
\PicA(\A) \cong \Pic(M) \cong H^2(M,\mathbb{Z}),
\end{equation}
where $\Pic(M)$ is the set of isomorphism classes of complex line bundles over $M$,
with group operation given by tensor product.
The first isomorphism in \eqref{eq:identif} is a consequence of the
Serre-Swan theorem (see e.g. \cite{Bass68}), 
and the second is the Chern class map.
Since all the automorphisms of $C^\infty(M)$ come from 
diffeomorphisms\footnote{We do not know of any published source for this ``folk
  theorem''.  The key step is to show that any homomorphism from $C^\infty(M)$
to $\complex$ is given by evaluation at a point.  The complex case
  which interests us here can be reduced to the real case by a 
  simple, elementary  argument (shown to us by Bill Arveson).  
  The case of real-valued functions is standard, though Maksim Maydanskiy
  has pointed out to us that the correspondence may fail if the set of
  components of $M$ has cardinality greater than that of the continuum.
  See Chapter IV of \cite{KrMi} for an extensive discussion of these
  issues.} of $M$,   
\eqref{eq:exactsq2} becomes for this example
\begin{equation}\label{eq:exactH2alg}
1 \rightarrow H^2(M,\mathbb{Z}) \rightarrow \Pic(C^\infty(M)) \stackrel{\proj}{\rightarrow} \Diff(M),
\end{equation}
and we obtain the purely geometric description of $\Pic(C^\infty(M))$
as the semidirect product $\Diff(M) \ltimes \Pic(M)$, with action given by
pull-back of line bundles.
\end{example}

In addition to the homomorphisms \eqref{eq:hhom}, one has
for any invertible morphism between $\A$ and $\B$ in $\catalg$ 
an algebra isomorphism between $\cent(\A)$ and $\cent(\B)$.
In particular,  two unital {\it commutative} $k$-algebras 
are isomorphic in $\catalg$ (i.e., Morita equivalent)
if and only if they are isomorphic in the usual sense \cite{Bass68}.
Despite this, regarding commutative
algebras as objects in $\catalg$ 
has the effect of enlarging the possible ways an algebra can be
``isomorphic'' to itself. 
As we saw in Example \ref{ex:linebundles}, for $\A=C^\infty(M)$
 these ``extra ways'' (i.e., $\PicA(\A)$) 
have a geometric interpretation as
the (isomorphism classes of) line bundles over $M$.

As we will discuss in the next sections, groupoids, symplectic
groupoids, and integrable Poisson manifolds can be similarly regarded as objects in
more general categories, in which the morphisms are geometric ``bimodules''.
Inspired by Example \ref{ex:linebundles}, 
we shall investigate Picard groups in this setting, the 
contribution to them of the geometric automorphisms, and
the  analogues of the exact sequences (\ref{eq:exactsq}) and \eqref{eq:exactsq2}.
 
\section{Morita equivalence and Picard groups of Lie groupoids}

\subsection{Generalized morphisms of Lie groupoids}\label{sec:genmorgr}
Generalized morphisms of Lie groupoids were studied in detail by
Mr{\v{c}}un \cite{Mrcun99},  who
called them {\bf Hilsum-Skandalis maps}, following \cite{HS87} 
(see also \cite{MM03} and references therein).  We
begin our discussion by fixing some notation and terminology.

Given a Lie groupoid $\Gamma$ over a
manifold $P$, we denote its
unit map by $\varepsilon: P \rightarrow \Gamma$, the inversion by
$i: \Gamma \rightarrow \Gamma$ and the target (resp.~source) map
by $t: \Gamma \rightarrow P$ (resp.~$s:\Gamma \rightarrow P$).
We will often identify an 
element $x$ in $P$ with its image $\varepsilon(x) \in \Gamma$.
 A {\bf (left) action}
of $\Gamma$ on a manifold $S$ consists of a map $P\stackrel{J}{\leftarrow} S$ called the
{\bf moment} and a map from $\{(g,x) \in \Gamma\times S \;|\;
s(g)=J(x)\}$ to $S$ satisfying the appropriate axioms (see for example
\cite{MiWe}).  Right actions are defined in a similar way (in which
case we indicate the moment as $S\stackrel{J}{\rightarrow} P$).  The
action is {\bf principal} with respect to a mapping $p:S\to M$ if 
$p$ is a surjective submersion and if $\Gamma$ acts freely and transitively on
each fibre of $p$.  ($p:S\to M$ is sometimes referred to as a
principal $\Gamma$-bundle.)

If the groupoids $\Gamma_1$ over $P_1$ and $\Gamma_2$ over $P_2$ act
on a manifold $S$ from the left and right, respectively, and the
actions commute,
we call $S$ a $(\Gamma_1,\Gamma_2)$-{\bf bibundle}.  Let $\PSP$
be the moments.
The bibundle
induces a correspondence between the orbit spaces of $\Gamma_1$ and
$\Gamma_2$ by the rule that orbits $\calo_1$ and $\calo_2$ are related if
$J_1^{-1}(\calo_1)$ meets $J_2^{-1}(\calo_2).$  

Since the essential ``ingredients'' of a groupoid are its orbits and its
isotropy groups,
a bibundle should represent a
  generalized morphism $\Gamma_1\leftarrow\Gamma_2$ only when 
the relation induced by $S$ is a map of orbit spaces
  $P_1/\Gamma_1\leftarrow P_2/\Gamma_2,$  and when each point
$x\in S$ determines a homomorphism from the isotropy group of
$J_2(x)$ in $\Gamma_2$ to that of $J_1(x)$ in $\Gamma_1$.  Both of
these conditions are met when
   the bibundle is {\bf left principal}, i.e., when the left action of
  $\Gamma_1$ is principal with respect to the right moment map $J_2$.
  (Transitivity of the $\Gamma_1$ action on $J_2$ fibres
makes the map on orbit spaces single-valued, while freeness gives the action on 
isotropy groups.)
It follows from the assumption of left principality 
(see \cite{Land01}, Lemma 4.18) that the
$\Gamma_1$ action is proper and that the map between orbit spaces is
smooth when the orbit spaces are manifolds.  

It is sometimes useful to think of each element $x$ of $S$ as an arrow
pointing from $J_2(x)$ to $J_1(x)$.  The operations of $\Gamma_1$ and
$\Gamma_2$ on $S$ may then be thought of as composition with the
arrows in the groupoids.  $S$ is left principal if the ``source'' map
$J_2$ is a surjective submersion, and if any two arrows with the same
source have a unique ``quotient'' in $\Gamma_1$.  
In other words, the bimodule is essentially the same thing as 
a (small) category whose set of objects is the
disjoint union of $P_1$ and $P_2$, whose restriction to each $P_i$ is
the groupoid $\Gamma_i$, and whose morphisms from $P_2$ to $P_1$
(there are none in the other direction) have a certain divisibility
property.

A ``tensor product'' of  bibundles \cite{HS87}
is then defined by ``composition of arrows'' as follows.  Let $\Gamma_j$
be a groupoid over $P_j$, for $j=1,2,3.$
If $S$ is a 
$(\Gamma_1,\Gamma_2)$-bibundle with moments $\PSP$, and $S'$ is
a $(\Gamma_2,\Gamma_3)$-bibundle with moments $\PSPprime$, their 
product is the orbit space
$$
S*S' := (S \times_{\st{P_2}} S')/{\Gamma_2},
$$
where $S \times_{\st{P_2}} S' = \{ (x,y) \in S \times 
S' \; | \; J_2(x)=J_2'(y)\}$
and the (right) $\Gamma_2$-action is given by 
$$
(x,y)\stackrel{g}{\mapsto} (xg,g^{-1}y).
$$
Fortunately, the assumption that $S$ and $S'$ be
left principal insures that $S*S'$ is a smooth manifold,
and that
$$
 \xymatrix {
{\Gamma_1} \ar@<-0.5ex>[d]\ar@<0.5ex>[d] & S*S'\ar[dl]^{J_1}\ar[dr]_{J_3'} &
{\Gamma_2} \ar@<-0.5ex>[d]\ar@<0.5ex>[d]\\
P_1 & & P_3
 }
$$ 
is a left principal $(\Gamma_1,\Gamma_3)$-bibundle.

Two $(\Gamma_1,\Gamma_2)$-bibundles $S_1$ and $S_2$ are {\bf isomorphic} if
there is a diffeomorphism between them which commutes with the
groupoid actions (including their moments).  It turns out that the
product $*$ is associative up to (natural) isomorphism, so that we may
define a category $\catgr$ in which the objects are Lie groupoids
and the morphisms $\Gamma_1\leftarrow\Gamma_2$ are isomorphism classes
of left principal $(\Gamma_1,\Gamma_2)$-bibundles, which we call {\bf generalized
  morphisms} between groupoids.  Following the terminology for algebras, 
we call two Lie groupoids $\Gamma_1$ and $\Gamma_2$
{\bf Morita equivalent} if they are isomorphic as objects in $\catgr$; 
this condition is equivalent to the existence of a $(\Gamma_1,\Gamma_2)$-bibundle
which is left and right principal (i.e., {\bf biprincipal}). Such a bibundle
is also called a {\bf Morita bibundle}.

\begin{example}(Gauge groupoids)\label{ex:gaugegrp}

Let $\Gamma$ be a Lie groupoid over $P$, and let $E \stackrel{p}{\to}
B$ be a right principal $\Gamma$-bundle.  The {\bf gauge groupoid}
associated with this principal bundle is the quotient of the pair
groupoid $E \times E$ by the right diagonal action of $\Gamma$. More
explicitly, it is a groupoid over $B$, with target (resp.~source) map
given by $[(x,y)] \mapsto p(x)$ (resp.~$[(x,y)] \mapsto p(y)$), and
multiplication defined by $[(x_1,y_1)][(x_2,y_2)]=[(x_1,y_2)]$, where
we choose representatives so that $y_1 = x_2$ (which is always
possible since $\Gamma$ acts transitively on $p$-fibres). We denote
the gauge groupoid of $E$ by $\mathcal{G}(E)$.

Let us assume that the moment $J$ for the $\Gamma$-action on $E$ is a surjective
submersion.
Note that $E$ carries a natural left $\mathcal{G}(E)$-action, with moment $p:E \to B$, defined by
$[(x,y)]\cdot z = x$, where $x$ is uniquely determined by the condition $y=z$. A simple computation
shows that this action is principal with respect to $J:E \to P$, and it is easy to see
that it commutes with the right $\Gamma$-action. As a result, 
\begin{equation}\label{eq:gaugebibundle}
 \xymatrix {
{\mathcal{G}(E)} \ar@<-0.5ex>[d]\ar@<0.5ex>[d] & E\ar[dl]_{p}\ar[dr]^{J} &
{\Gamma} \ar@<-0.5ex>[d]\ar@<0.5ex>[d]\\
B & & P
 }
\end{equation}
is a biprincipal $(\mathcal{G}(E),\Gamma)$-bibundle, and $\mathcal{G}(E)$ and $\Gamma$
are Morita equivalent.
\end{example}

\subsection{Picard groups of Lie groupoids}
\label{subsec-piclie}
Let us recall that a groupoid 
homomorphism is a functor, when groupoids are thought of as categories.
Any groupoid homomorphism
$\Gamma_1 \stackrel{\Phi}{\leftarrow}\Gamma_2$
can be seen as a generalized morphism as follows: we let $(\Gamma_1)_{\Phi}$ consist of pairs
$(g,y) \in \Gamma_1 \times P_2$ such that $s(g)=\Phi(y)$; this space carries a left principal
bibundle structure defined by a left $\Gamma_1$-action, with moment $(g,y) \mapsto t(g)$,
by $g_1 \cdot (g,y) = (g_1 g,y)$, and a right $\Gamma_2$-action, with moment 
$(g,y) \mapsto y$, by $(g,y)\cdot g_2 = (g\Phi(g_2), s(g_2))$. 
One can check that this correspondence of morphisms preserves their composition,
\begin{equation}\label{eq:compos}
(\Gamma_1)_{\Phi}*(\Gamma_2)_{\Psi}\cong (\Gamma_1)_{\Phi\circ \Psi},
\end{equation} 
and gives rise to a functor $\inc$ from the ``conventional'' category of groupoids
and groupoid homomorphisms into $\catgr$.  (See, for example,
\cite{MM03}.  The case of \'etale groupoids is treated in
\cite{Mrcun99}, but the general case is essentially the same.)

Let $\Gamma$ be a Lie groupoid over $P$. The group of ordinary groupoid automorphisms of $\Gamma$ 
is denoted by $\Aut(\Gamma)$, and we define the {\bf Picard group} of $\Gamma$, $\Pic(\Gamma)$, as
the group of automorphisms of $\Gamma$ regarded as an object in $\catgr$ (i.e., $\Pic(\Gamma)$
is the group of isomorphism classes of biprincipal $(\Gamma,\Gamma)$-bibundles). The unit in $\Pic(\Gamma)$
is the isomorphism class of $\Gamma$, regarded as a $(\Gamma,\Gamma)$-bibundle with respect to left
and right multiplication.

As in the case of algebras, the functor $\inc$ restricts to a group
homomorphism
\begin{equation}\label{eq:jhomom}
\inc: \Aut(\Gamma) \longrightarrow \Pic(\Gamma), \;\; \Phi \mapsto [(\Gamma)_{\Phi}]. 
\end{equation}
We now discuss the analogue of the exact sequence \eqref{eq:exactsq}.

Recall that a submanifold $N$ of a groupoid $\Gamma$ over $P$
is called a {\bf bisection} if
the restrictions $t_N$ and $s_N$ of $t$ and $s$ are diffeomorphisms
from $N$ to $P$.
The bisections form a group $\Bis(\Gamma)$ with the operations of
setwise multiplication and
inversion.  A bisection $N$
induces diffeomorphisms (called left
and right {\bf slidings} by Albert and Dazord \cite{AlDa91}) $l_N$ and $r_N$
of $\Gamma$ by  
$l_N(g)=ag$ where $a \in N$ is such 
that $s(a)=t(g)$, and $r_N(g)=gb$ where $b\in N$ is such that
$s(g)=t(b)$.  Left slidings commute with right slidings, and the
composition of $l_N$ with $r_{N^{-1}}$ is an automorphism $\Phi_N$,
  which we call an {\bf inner automorphism} of $\Gamma$. 
(Here $N^{-1}=i(N)$, where $i$ is the inversion map.)
The group of inner automorphisms of $\Gamma$ is denoted by $\Inaut(\Gamma)$.

\begin{proposition}\label{prop:exactsqgr}
Let $\Gamma$ be a Lie groupoid over $P$. The following is an exact sequence:
\begin{equation}\label{eq:exactsqgr}
1 \to \Inaut(\Gamma) \to \Aut(\Gamma) \stackrel{\inc}{\to} \Pic(\Gamma).
\end{equation}
\end{proposition}
\begin{proof}
Let $\Phi \in \Aut(\Gamma)$.
Note that $(\Gamma)_{\Phi}$ can be identified with the manifold $\Gamma$,
carrying a left $\Gamma$-action by left multiplication and a right
$\Gamma$-action, with moment $\Phi^{-1}\circ s : \Gamma \to P$, by
$z \cdot g = z\Phi(g)$.

Suppose $f: (\Gamma)_{\Phi} \to \Gamma$ is a $(\Gamma,\Gamma)$-bibundle isomorphism.
Since $f$ commutes with moments, we must have $t(f(x))=x$ and $s(f(x))=\Phi^{-1}(x)$,
for all $x \in P$. It then follows that the submanifold $N=f(P)$ is a bisection
of $\Gamma$.

Now, if $z \in \Gamma$ and $x=s(z)$, then $f(z)=f(z x)= z a$, where $a = f(x) \in N$
is uniquely determined by $t(a)=s(z)$. Thus $f = r_N$.  On the other hand, since $f$
is a bibundle isomorphism, we have $f(z\Phi(g)) = z\Phi(g)b$, for $b \in N$ satisfying
$t(b)=s(\Phi(g))$, and also $f(z\Phi(g))=f(z)g = zb'g$, for $b' \in N$ such that
$t(b')=s(z)=t(\Phi(g))$. Note that $b'$ satisfies $s(b')=t(g)$, and this condition
also determines it uniquely. So $\Phi(g)b = b'g$, which is equivalent to 
$\Phi = l_N\circ {r_N}^{-1} = \Phi_N$.

Conversely, if $\Phi_N$ is an inner automorphism of $\Gamma$, it is easy to check that
the right sliding $r_{N}$ is a bibundle isomorphism $(\Gamma)_{\Phi} \to \Gamma$.
\end{proof}

Let $\CIsoBis(\Gamma)$ be the subgroup of bisections with values in the centers of
the isotropy subgroups of $\Gamma$. This is the subgroup of $\Bis(\Gamma)$
giving rise to trivial inner automorphisms; i.e., 
the following is an exact sequence:
\begin{equation}\label{eq:exactcenter}
1 \rightarrow \CIsoBis(\Gamma) \rightarrow \Bis(\Gamma) \rightarrow 
\Inaut(\Gamma) \rightarrow 1.
\end{equation}

Let $\Outaut(\Gamma):= \Aut(\Gamma)/\Inaut(\Gamma)$ denote the group of {\bf outer automorphisms} of
$\Gamma$. The next simple observation describes when $\Pic(\Gamma)$ coincides
with $\Outaut(\Gamma)$.

\begin{lemma}\label{lem:jonto}
The exact sequence \eqref{eq:exactsqgr} extends to
\begin{equation}\label{eq:jonto}
1 \to \Inaut(\Gamma) \to \Aut(\Gamma) \stackrel{\inc}{\to} \Pic(\Gamma) \to 1
\end{equation}
(i.e., $\inc$ is onto) if and only if for  every biprincipal $(\Gamma,\Gamma)$-bibundle
$S$, with moments $\PSPn$, there exists a smooth map
$\sigma:P \to S$ with $J_2\circ \sigma =\id$ and $J_1\circ \sigma \in \Diff(P)$. 
In this case, $\Pic(\Gamma)\cong \Outaut(\Gamma)$.
\end{lemma}

\begin{proof}
Let $S = (\Gamma)_{\Phi}$, for $\Phi \in \Aut(\Gamma)$. Then $\sigma=\varepsilon \circ \Phi|_P$, where
$\varepsilon$ is the identity embedding, is a cross section for the moment 
$J_2=\Phi^{-1}\circ s$, and $J_1 \circ \sigma =
t \circ \varepsilon \circ \Phi|_P = \Phi|_P  \in \Diff(P)$. 

Conversely, let $S$ be a $(\Gamma,\Gamma)$-Morita bibundle, with moments $J_1$
and $J_2$, and suppose $\sigma:P \to S$ is a map satisfying 
$J_2\circ \sigma =\id$ and $J_1\circ \sigma \in \Diff(P)$. Let us define the map
$f:S \to \Gamma$ by
$f(z) = g$, where $g$ is the unique element in $\Gamma$ satisfying
$g \sigma(J_2(z)) = z$ in $S$. It is clear that $f$ preserves left $\Gamma$-actions,
and it is a diffeomorphism as a result of the principality of the left
$\Gamma$-action on $S$ and the condition $J_1\circ \sigma \in \Diff(P)$.
Note that the condition $f(z\cdot g) = f(z)\Phi(g)$ defines an
automorphism $\Phi\in \Aut(\Gamma)$ so that $S \cong (\Gamma)_{\Phi}$
as a bibundle. So $\inc$ is onto.
\end{proof}

Let us give some concrete examples.

\begin{example}(Lie groups)\label{ex:liegrp}

Let $G$ be a Lie group, regarded as a Lie groupoid over a point.  ($G$
need not be connected; in particular, it may be discrete.)
Let $\X$ be a biprincipal $(G,G)$-bibundle, also called a $G$-{\bf bitorsor},
defining a Morita self-equivalence of $G$. 
A map $\sigma$ from the base of $G$ into $\X$ is just a choice of a point in $\X$.
It immediately follows from  Lemma \ref{lem:jonto} that there exists
$\phi \in \Aut(G)$ such that
$\X$ is isomorphic, as a bibundle, to $G_{\phi}$ (action on the left
by left multiplication, and on the right by right multiplication composed with $\phi$).
Hence the exact sequence \eqref{eq:exactsqgr} becomes
$$
1 \to \Inaut(G) \to \Aut(G) \stackrel{\inc}{\to} \Pic(G) \to 1,
$$
so $\Pic(G)\cong \Outaut(G)$.  
\end{example}

Since Morita equivalence is just isomorphism in the category $\catgr$, it is clear that
Morita equivalent groupoids have isomorphic Picard groups.
As a consequence of Examples \ref{ex:gaugegrp} and \ref{ex:liegrp}, we obtain

\begin{corollary}\label{cor:moritaexample}
If $\mathcal{G}(E)$ is the gauge groupoid of a principal $\Gamma$-bundle $E$ (so that
its moment $J$ is a surjective submersion), then
$\Pic(\mathcal{G}(E)) \cong \Pic(\Gamma)$. In particular, if $\Gamma=G$
is a Lie group, then $\Pic(\mathcal{G}(E)) \cong \Outaut(G)$.
\end{corollary}

\begin{example}\label{ex:gendeg}(Transitive groupoids)

A transitive groupoid $\Gamma$ over $P$ is always
isomorphic to a gauge groupoid of a Lie group principal bundle:
for a fixed point $x \in P$, we consider
the  set of arrows starting at $x$, $E_x := s^{-1}(x)$;
if $\Gamma_x$ denotes the isotropy subgroup of $\Gamma$ at $x$,
a simple computation shows that $E_x$ is a principal
$\Gamma_x$-bundle, in such a way that $\Gamma$ is isomorphic
to the gauge groupoid $\mathcal{G}(E_x)$. 
It follows from Corollary \eqref{cor:moritaexample} that 
$$
\Pic(\Gamma)\cong \Pic(\mathcal{G}(E_x))\cong \Outaut(\Gamma_x).
$$
In particular, this shows that, although the isotropy groups of
$\Gamma$ at different points of $P$ 
are isomorphic to one another only in a noncanonical way,
their groups of outer isomorphisms can be canonically identified with
each other.  

The homomorphism $\inc:\Aut(\Gamma)\to \Outaut(\Gamma_x)$ can be described directly
as follows.
An automorphism $\Phi \in \Aut(\Gamma)$
induces an isomorphism $\Phi_x: \Gamma_x \to \Gamma_{\Phi(x)}$. Since $\Gamma$
is transitive, one can identify $\Gamma_x$ and $\Gamma_{\Phi(x)}$ by means
of a choice of $g$ with $s(g)=x$, $t(g)=\Phi(x)$. This induces an
automorphism 
$$\Phi_x^g:\Gamma_x \to \Gamma_x, \;\;\; \Phi_x^g(z)=g^{-1}\Phi_x(z)g,
$$ and a simple
computation shows that, if $h \in \Gamma$ also satisfies $s(h)=x$, $t(h)=\Phi(x)$,
then $\Phi_x^g\circ (\Phi_x^h)^{-1} \in \Inaut(\Gamma_x)$. So $\Phi$ determines an
equivalence class $[\Phi_x] \in \Outaut(\Gamma_x)$, which is the image
of $\Phi$ under $j$.
So, for a transitive groupoid, the exact sequence \eqref{eq:exactsqgr} becomes
\begin{equation}
1 \to \{\Phi \in \Aut(\Gamma) \,|\, [\Phi_x]=\id \in \Outaut(\Gamma_x)\}
\to \Aut(\Gamma) \to \Outaut(\Gamma_x).
\end{equation}
In this (indirect) way, we get a characterization of inner automorphisms of transitive
groupoids as
\begin{equation}\label{eq:inndeg}
\Inaut(\Gamma)
= \{\Phi \in \Aut(\Gamma) 
\,|\, [\Phi_x]=\id \in \Outaut(\Gamma_x)\}.
\end{equation}
This identification  can also be verified directly:
if $\Phi_N$ is an inner automorphism associated with a bisection $N$
and  we consider $(\Phi_N)_x:\mathcal{G}(E)_x \to 
\mathcal{G}(E)_{\Phi_N(x)}$, $x \in P$,
then $(\Phi_N)_x^g$ is trivial if we choose $g \in N$; conversely,
if $[\Phi_x]=0$, we can reconstruct $N$ (uniquely up to $\CIsoBis(\Gamma)$, see
\eqref{eq:exactcenter}) as elements $g$ in $\Gamma$ so that $(\Phi_N)_x^g$ is trivial.
\end{example}

\begin{example}\label{ex:funda}(Fundamental groupoids)

Let $P$ be a connected manifold.  The {\bf fundamental groupoid}
$\Pi(P)$, consisting of homotopy classes of paths in $P$ with fixed
endpoints, is a transitive groupoid over $P$ whose isotropy group at
$x\in P$ is the fundamental group $\pi_1(P,x).$  The corresponding
principal bundle is the universal covering $\widetilde{P}\rightarrow
P$.  Thus, the outer automorphism groups $\Outaut(\pi_1(P,x))$ are all
canonically isomorphic to $\Pic(\Pi(P))$ and hence to one another.
\end{example}

We now discuss the groupoid analogue of \eqref{eq:exactsq2}.
For a groupoid $\Gamma$ over $P$, the orbit space $\cent(\Gamma)=P/\Gamma$ will
play the role of its ``center,'' since the inner automorphisms act
trivially on it.   In this sense, a ``commutative'' groupoid
is a groupoid with trivial orbits ($\cent(\Gamma) = P$), and the groupoids
with trivial center are the transitive ones.

In general, we regard the orbit space $\cent(\Gamma)$ (the ``center''
of $\Gamma$) as a topological space, with the quotient
topology. However, whenever $\cent(\Gamma)$ is smooth and such that
the quotient map is a submersion, we regard it as a smooth manifold.
A $(\Gamma_1,\Gamma_2)$-Morita bibundle induces a homeomorphism (or
diffeomorphism, in the smooth case) between $\cent(\Gamma_1)$ and
$\cent(\Gamma_2)$, and therefore defines a homomorphism
\begin{equation}\label{eq:hhomgr}
h: \Pic(\Gamma) \rightarrow \Aut(\cent(\Gamma)).
\end{equation}
Here $\Aut(\cent(\Gamma))$ denotes the group of homeomorphisms 
(or diffeomorphisms, in the smooth case) of $\cent(\Gamma)$ onto itself. 
The kernel of $h$ consists of those bibundles which induce 
the identity map on $\cent(\Gamma)$;
we denote this subgroup of $\Pic(\Gamma)$ by $\PicZG(\Gamma)$.  It
seems appropriate to
call this the {\bf static Picard group} of $\Gamma$. 
The analogue of \eqref{eq:exactsq2} is
\begin{equation}\label{eq:exactgr2}
1 \rightarrow \PicZG(\Gamma) \rightarrow \Pic(\Gamma) \stackrel{\proj}{\rightarrow} \Aut(\cent(\Gamma)).
\end{equation}

The next two examples discuss the Picard groups of Lie groupoids $\Gamma$
over $P$ with $\cent(\Gamma)=P$ (i.e. $s=t$). These groupoids are
just bundles of Lie groups, i.e.,  smooth families of Lie groups parametrized by $P$.
The description of the Picard groups of such groupoids 
follows from the study of their bitorsors 
by Moerdijk in \cite{Moerd02}.

\begin{example}\label{ex:trivcontr}(Bundles of contractible Lie groups)

Let $\Gamma$ be a bundle of contractible Lie groups
over $P$.  
Let $S$ be a $(\Gamma,\Gamma)$-Morita bibundle, with
moments $J_1$ and $J_2$. The $J_2$-fibres are contractible (since each $J_2$-fibre
is  diffeomorphic to an $s$-fibre), so there exists a cross-section $\sigma:P\to S$,
$J_2\circ \sigma = \id$. Note that $J_1\circ \sigma$ coincides with the diffeomorphism
induced by $S$ on $\cent(\Gamma)=P$.
So, by Lemma \ref{lem:jonto}, 
$\Pic(\Gamma) \cong \Outaut(\Gamma)$, which consists of automorphisms
of the family of groups, modulo smooth sections acting by inner
automorphisms.  
The static Picard group is thus
\begin{equation}\label{eq:picztriv}
\PicZG(\Gamma)\cong \{[\Phi] \in \Outaut(\Gamma) \;|\; \Phi|_{P}=\id\}.
\end{equation}

In particular, if $\Gamma$ is a bundle of contractible {\it abelian} Lie groups 
(e.g. a vector bundle with groupoid
structure given by fibrewise addition), inner automorphisms are
trivial, and so $\Pic(\Gamma) \cong \Aut(\Gamma)$.
\end{example}

\begin{example}\label{ex:abelian}(Bundles of abelian Lie groups)

This example follows \cite[Sec.~3.5]{Moerd02}.

Let $\Gamma$ be a bundle of abelian Lie groups.
In this case, one can show that any $(\Gamma,\Gamma)$-Morita
bibundle is completely determined by a pair $(S,\Phi)$, where $S$ is
a principal $\Gamma$-bundle
(also called a $\Gamma$-{\bf torsor})
and $\Phi \in \Aut(\Gamma)$. Moreover,
there is a natural tensor product operation on $\Gamma$-torsors making
their set of isomorphism classes into a group isomorphic to
$\check{H}^1(P,\Gamma)$. The Picard group of $\Gamma$ can then
be identified with the semi-direct product
$$
\Pic(\Gamma)\cong \Aut(\Gamma)\ltimes \check{H}^1(P,\Gamma).
$$
It also follows that 
$$
\PicZG(\Gamma) \cong \Aut_P(\Gamma)\ltimes \check{H}^1(P,\Gamma),
$$
where $\Aut_P(\Gamma)=\{\Phi \in \Aut(\Gamma)\,|\, \Phi|_{P}=\id\}$.
\end{example}

\section{Symplectic dual pairs as generalized morphisms}

Let $S$ be a symplectic manifold. A pair of Poisson maps
$$ 
\begin{diagram}
\node{}\node{S}\arrow{sw,t}{J_1}\arrow{se,t}{J_2}\node{}\\
\node{P_1}\node{}\node{P_2}
\end{diagram}
$$
is called a {\bf dual pair} if the $J_1$- and $J_2$-fibres are the symplectic
orthogonal of one another. 
A dual pair is called {\bf full} if $J_1$ and $J_2$ are surjective submersions;
it is called {\bf complete}
if the maps $J_1$ and $J_2$ are complete (recall that a Poisson map
$J:Q \to P$ is {\bf complete} if the hamiltonian vector field
$X_{J^*f}$ is complete whenever $X_f$ is complete, $f \in C^\infty(M)$).

We now discuss how dual pairs relate to generalized morphisms
of symplectic groupoids and Poisson manifolds.

A {\bf symplectic groupoid} \cite{CDW87,We87} is a Lie groupoid
$\Gamma$ equipped with a symplectic structure $\omega$ for which the graph
of the multiplication $m: \Gamma_2 \rightarrow \Gamma$ is lagrangian
in $\Gamma \times \Gamma \times \overline{\Gamma}$ (where
$\overline{\Gamma}$ is equipped with $-\omega$). 
This compatibility condition between the groupoid structure
and the symplectic form implies, in particular, that the identity embedding
$\varepsilon:P \to \Gamma$ is lagrangian, the inversion $i:\Gamma \to \Gamma$
is an anti-symplectomorphism and $s$- and $t$-fibres are symplectically
orthogonal to one another; furthermore, 
the identity section $P$ inherits a Poisson structure $\pi$, uniquely
determined by the condition that $t:\Gamma \to P$ (resp.~$s:\Gamma \to P$)
is a Poisson map (resp.~anti-Poisson map).

Let $(S,\omega_S)$ be a symplectic manifold. A (left) action of $(\Gamma,\omega)$
on $S$ is called {\bf symplectic} if the graph of the action map
$\Gamma \times S \to S$ is lagrangian in $\Gamma \times S \times \overline{S}$.
In this case, the moment $J:S \to P$ is automatically
a Poisson map \cite{MiWe}. Symplectic right actions are defined
analogously, but their moments are anti-Poisson maps.  

If $\Gamma_1$ and $\Gamma_2$ are symplectic groupoids,
then a $(\Gamma_1,\Gamma_2)$-bibundle $S$ is a {\bf symplectic bibundle}
if both actions are symplectic. 
An isomorphism of symplectic bibundles
is just an isomorphism of bibundles preserving the symplectic structures.
A {\bf generalized morphism} 
$\Gamma_1 \leftarrow \Gamma_2$ is an
isomorphism class of
symplectic $(\Gamma_1,\Gamma_2)$-bibundles which are left principal. 
It is proven in \cite{Xu91b} that the ``tensor product''
of symplectic bibundles (regarded just as bibundles, see Section \ref{sec:genmorgr}) 
is automatically compatible
with the symplectic structures, in such a way that the resulting
bibundle is canonically symplectic.
Thus, the tensor product operation induces a well-defined composition
of generalized morphisms; we denote the resulting category by
$\catsym$.

We call two symplectic groupoids $\Gamma_1$ and $\Gamma_2$ {\bf Morita equivalent}
if they are isomorphic in $\catsym$, or, equivalently, if there exists a biprincipal
symplectic $(\Gamma_1,\Gamma_2)$-bibundle \cite{Land00b}. One can check that,
in this case, the moments
$\PSPb$ define a complete and full dual pair (here $\overline{P}_2$ denotes
the manifold $P_2$ with the opposite Poisson structure).

As mentioned earlier in this section, the identity section $P$
of a symplectic groupoid $(\Gamma,\omega)$ naturally inherits
a Poisson structure making $t$ into a Poisson map (and $s$ into
an anti-Poisson map).
A Poisson manifold $(P,\pi)$
is called {\bf integrable} if there exists a symplectic groupoid
over $P$ for which the induced Poisson structure on $P$ is $\pi$.
When $(P,\pi)$ is integrable, it has a canonically defined
source-simply-connected\footnote{By ``simply-connected,'' we will always
  mean connected, with trivial fundamental group.} symplectic groupoid
$\Gamma(P)$.  
Not every Poisson manifold is integrable, and the obstructions have been
explicitly described in \cite{CrFe02} (we will return to this topic in
Section \ref{sec:extended}). A Poisson structure $\pi$ on $P$
defines a Lie algebroid structure on $T^*P$ (see e.g. \cite{CDW87}),
and $P$ is integrable if and only if $T^*P$ is integrable
as a Lie algebroid \cite{MaXu} (see also \cite{CrFe02}).

Symplectic actions of symplectic groupoids and symplectic bibundles
can be described purely in terms of the Poisson geometry of the
identity sections.  Recall that a {\bf symplectic realization} 
of a Poisson manifold $P$ is just a Poisson map from a symplectic
manifold to $P$.  
Any symplectic realization $J: S \rightarrow P$ induces a canonical
action of the Lie algebroid $T^*P$ (induced by $\pi$) on $S$ by
assigning to each $1$-form $\alpha$ on $P$ the vector field $X$ on $S$
defined by $i_X\omega_{\st{S}} = J^*\alpha$. If $P$ is integrable,
this action
extends to a symplectic action of $\Gamma(P)$ when $J$ is complete
\cite{CDW87,CrFe02}, in which case $J$ is the moment of the action.  
Hence there is a natural correspondence between complete
symplectic realizations of $P$ and symplectic actions of $\Gamma(P)$.

If $P_1$ and $P_2$ are integrable Poisson manifolds and
\begin{equation}\label{eq:symbim}
\begin{diagram}
\node{}\node{S}\arrow{sw,t}{J_1}\arrow{se,t}{J_2}\node{}\\
\node{P_1}\node{}\node{\overline{P_2}}
\end{diagram}
\end{equation}
is a pair of complete symplectic realizations,
then $S$ carries a left $\Gamma(P_1)$-action and a right $\Gamma(P_2)$-action.
These actions commute if and only if they commute on the infinitesimal
level, i.e.,
\begin{equation}\label{eq:commute}
\{J_1^*C^\infty(P_1), J_2^*C^\infty(P_2)\} =0.
\end{equation}
In this case, we call \eqref{eq:symbim} a $(P_1,P_2)$-{\bf symplectic bimodule}.
It the follows that there is a natural
one-to-one correspondence between $(P_1,P_2)$-symplectic bimodules and 
$(\Gamma(P_1),\Gamma(P_2))$-symplectic bibundles.
A {\bf generalized morphism} of integrable Poisson manifolds $P_1$ and $P_2$
is an isomorphism class of $(P_1,P_2)$-symplectic bimodules
which correspond to left principal $(\Gamma(P_1),\Gamma(P_2))$-symplectic
bibundles. 

Following the ideas in \cite[Thm.~3.2]{Xu91}, one obtains 
an alternative way to define generalized
morphisms of Poisson manifolds just in terms of Poisson structures and moments,
with no reference to symplectic
groupoids: a symplectic bimodule $\PSP$ is a generalized morphism if and only if
$J_1$ and $J_2$ are complete Poisson maps, $J_1$ is a submersion, $J_2$ is a surjective
submersion with simply-connected fibres, and the $J_1$- and $J_2$-fibres are
the symplectic orthogonals of each other.
The category whose objects are integrable Poisson manifolds
and arrows are generalized morphisms is denoted by $\catpoiss$.

\begin{remark}\label{rem:etale}
Any ordinary complete Poisson map
$P_1\stackrel{\phi}{\leftarrow} P_2$ defines a symplectic
bimodule, namely
$P_1\stackrel{\phi\circ  t_2}{\leftarrow}
\Gamma(P_2)\stackrel{s_2}{\rightarrow}P_2$; we denote this symplectic bimodule
by $_{\phi}(\Gamma(P_2))$.
It follows from our description of generalized morphisms of Poisson manifolds that
such a symplectic bimodule is a generalized morphism if and only if $\phi$
is \'etale, which means that $\phi$ is a submersion between leaf spaces and a covering 
on each symplectic leaf.

If $\phi$ is \'etale, it lifts to a Lie algebroid morphism $T^*P_1 \leftarrow
T^*P_2$ and thus to a symplectic groupoid morphism 
$\Gamma(P_1)\stackrel{\Phi}{\leftarrow} \Gamma(P_2)$. 
Hence $\phi$ also defines a generalized morphism $(\Gamma(P_1))_{\Phi}$, as discussed
in the beginning of Section \ref{subsec-piclie}. We note that there is a natural
$(P_1,P_2)$-bimodule map $_{\phi}(\Gamma(P_2)) \to (\Gamma(P_1))_{\Phi}$ given by
$g \mapsto (\Phi(g),s_2(g))$ which is \'etale;  it becomes 
an isomorphism whenever $\phi$ is a Poisson
diffeomorphism.
\end{remark}

We call two integrable Poisson manifolds
$P_1$ and $P_2$ {\bf Morita equivalent} \cite{Xu91} if their source-simply-connected symplectic
groupoids are Morita equivalent; this coincides with the notion of isomorphism
in $\catpoiss$ \cite{Land00b}. Equivalently, we can define Morita equivalence of Poisson manifolds
in Poisson geometrical terms \cite[Thm.~3.2]{Xu91}: $P_1$ and $P_2$ are Morita equivalent
if there is a symplectic manifold $S$ with complete Poisson and anti-Poisson maps
$J_1:S \rightarrow P_1$ and $J_2:S \rightarrow P_2$, so that $\PSPb$
is a complete full dual pair for which $J_1$- and $J_2$-fibres are
simply connected.

In this case, the diagram
$$ 
\begin{diagram}
\node{}\node{S}\arrow{sw,t}{J_1}\arrow{se,t}{J_2}\node{}\\
\node{P_1}\node{}\node{P_2}
\end{diagram}
$$
is called a {\bf Morita bimodule}.

Morita equivalent Poisson manifolds share many properties. For example, they
have isomorphic first Poisson cohomology groups \cite{GinzLu92}, homeomorphic
spaces of symplectic leaves \cite{BuRad02,CrFe02} and isomorphic 
transverse geometry at corresponding symplectic leaves \cite{We83}. 
For a symplectic manifold, its fundamental group is a complete
Morita invariant \cite{Xu91}; in particular, simply-connected symplectic manifolds are Morita equivalent
to a point. We will discuss further
examples of complete invariants of Morita equivalence in Section
\ref{sec:TSS1}.

\begin{remark}
\label{remark:restrict}
For later use in Section \ref{sec:TSS1}, we note that symplectic
leaves that correspond to each other in a Morita equivalence are
themselves Morita equivalent \cite{CrFe02}.  When the symplectic
leaves are open, one may simply restrict the equivalence.  (In fact,
one may also restrict a Morita equivalence to any open subset which is
a union of symplectic leaves.)  In general, consider symplectic leaves
$L_i$ in $P_i$, $i=1,2$, and $N= J_1^{-1}(L_1) = J_2^{-1}(L_2)$ (the
pull-back foliations by $J_1$ and $J_2$ coincide); then $J_i$ maps $N$
onto $L_i$, and the map has simply-connected fibers. So we have an
induced isomorphism of fundamental groups $\pi_1(L_1)\cong \pi_1(N)
\cong \pi_1(L_2).$
\end{remark}
\section{Picard groups of symplectic groupoids and Poisson manifolds}
\label{sec:picsymp}

Let $\Gamma$ be a symplectic groupoid. We denote its group of symplectic groupoid
automorphisms by $\Aut(\Gamma)$, 
and its group of automorphisms in $\catsym$
by $\Pic(\Gamma)$; similarly, if $P$ is an integrable Poisson manifold, the group of Poisson 
diffeomorphisms $f:P\to P$ is denoted 
by $\Poiss(P)$, whereas the group
of automorphisms of $P$ in $\catpoiss$ is denoted by $\Pic(P)$.
Note that, if $P$ is integrable, then we have a natural identification
\begin{equation}\label{eq:samepic}
\Pic(P)=\Pic(\Gamma(P)),
\end{equation}
where $\Gamma(P)$ is the canonical source-simply-connected symplectic groupoid integrating $P$.

As in the case of Lie groupoids \eqref{eq:jhomom}, we have a group homomorphism
\begin{equation}
j: \Aut(\Gamma) \to \Pic(\Gamma), \;\; \Phi \mapsto [(\Gamma)_{\Phi}],
\end{equation}
which induces, for Poisson manifolds, a group homomorphism
\begin{equation}\label{eq:poissmorph}
j: \Poiss(P) \to \Pic(P),\;\; \phi \mapsto [(\Gamma(P))_{\phi}],
\end{equation}
where $(\Gamma(P))_{\phi}$ is the symplectic bimodule
$$
\begin{diagram}
\node{}\node{\Gamma(P)}\arrow{sw,t}{t}\arrow{se,t}{\phi^{-1}\circ s}\node{}\\
\node{P}\node{}\node{P}
\end{diagram}
$$

On a symplectic groupoid $\Gamma$, 
the group of lagrangian bisections (i.e., bisections which are lagrangian
submanifolds), $\LBis(\Gamma)$, plays a special role: 
if $L$ is a lagrangian bisection, then
the inner automorphism
$\Phi_L$ is a symplectomorphism of $\Gamma$, and $\phi_L=(\Phi_L)|_{P} : P \to P$
preserves the induced Poisson structure on $P$ \cite{MiWe}.
We refer to the subgroup
$$
\Inaut(\Gamma):=\{ \Phi_L \,:\, L \in \LBis(\Gamma(P)) \} \subseteq \Aut(\Gamma)
$$
as the subgroup of {\bf inner automorphisms} of $\Gamma$; similarly,
the subgroup
$$
\InnPois(P):= \{ \phi_L=(\Phi_L)|_{P} \,:\, L \in \LBis(\Gamma(P)) \}\subseteq \Poiss(P)
$$
is called the subgroup of {\bf inner Poisson automorphisms} of $P$.
We denote the quotients  $\Poiss(\Gamma)/ \Inaut(\Gamma)$ by $\Outaut(\Gamma)$
and $\Poiss(P)/ \InnPois(P)$ by $\OutPois(P)$.  For example, when the
modular class of a Poisson manifold \cite{We97} is nonzero,
the flow of any modular vector field represents a nontrivial
one-parameter group of outer automorphisms.

The following result follows from Proposition \ref{prop:exactsqgr}.  

\begin{proposition}\label{prop:exactsqpoiss}
Let $\Gamma$ be a symplectic groupoid. Then the following is an exact sequence:
\begin{equation}\label{eq:exactpoiss0}
1 \to \Inaut(\Gamma) \to \Aut(\Gamma) \stackrel{j}{\to}\Pic(\Gamma).
\end{equation}
Similarly, if $P$ is a Poisson manifold, then the sequence below is exact:
\begin{equation}\label{eq:exactpoiss1}
1 \to \InnPois(P) \to \Poiss(P) \stackrel{j}{\to}\Pic(P).
\end{equation}
In particular, the outer automorphisms $\Outaut(\Gamma)$ (resp.~$\OutPois(P)$)
sit in $\Pic(\Gamma)$ (resp.~$\Pic(P)$) as a subgroup.
\end{proposition}

To describe the inner automorphisms in terms of the groupoid, we introduce
the subgroup $\IsoLBis(\Gamma(P))$ of $\LBis(\Gamma(P))$ consisting of 
bisections with values in the isotropy subgroupoid of $\Gamma(P)$, 
and its subgroup $\mathrm{ZIsoLBis}(\Gamma(P))$ consisting 
of lagrangian bisections with values in the centers of the isotropy subgroups.
We have the exact sequence corresponding to \eqref{eq:exactcenter}:
\begin{equation}\label{eq:exactbis}
1 \rightarrow \IsoLBis(\Gamma(P)) \rightarrow \LBis(\Gamma(P)) \rightarrow 
\InnPois(P) \rightarrow 1
\end{equation}
as well as a sequence which describes the kernel in
\eqref{eq:exactpoiss0}:
\begin{equation}\label{eq:exactbis2}
1 \rightarrow \mathrm{ZIsoLBis}(\Gamma(P)) \rightarrow \LBis(\Gamma(P)) \rightarrow 
\Inaut(\Gamma(P)) \rightarrow 1.
\end{equation}

If $\Gamma$ is a symplectic groupoid, then $\cent(\Gamma)$ is the space of
symplectic leaves of the Poisson structure on its identity section.
So, for an integrable Poisson manifold $P$, we define $\cent(P)$ as the leaf space $P/\mathcal{F}$,
where $\mathcal{F}$ is the symplectic foliation of $P$.
As in the case of Lie groupoids, we generally consider this leaf space as a topological space,
but we regard it as a smooth manifold whenever $\mathcal{F}$ is simple.
In this way, $\Aut(\cent(P))$ denotes the group of homeomorphisms (or diffeomorphisms, in the smooth
case) of $\cent(P)$ onto itself.
It follows from the discussion for Lie groupoids (see also \cite{BuRad02,CrFe02})
that a $(P_1,P_2)$-Morita bimodule
induces a homeomorphism (or diffeomorphism, in the smooth case)
between $\cent(P_1)$ and $\cent(P_2)$. As a result, we get a group homomorphism
\begin{equation}\label{eq:hhompois}
h: \Pic(P) \rightarrow \Aut(\cent(P)).
\end{equation}
The kernel of $h$ consists of those $(P,P)$-Morita  bimodules which induce 
the identity map on $\cent(P)$;
we denote this subgroup of $\Pic(P)$ by $\PicZP(P)$. We then obtain an analogue of \eqref{eq:exactsq2} and
\eqref{eq:exactgr2}:
\begin{equation}\label{eq:exactpois2}
1 \rightarrow \PicZP(P) \rightarrow \Pic(P) \stackrel{\proj}{\rightarrow} \Aut(\cent(P)).
\end{equation}

As with groupoids, we have two extreme situations: the ``commutative''
case, where $\cent(P) = P$, occurs for the zero Poisson
structure, while the ``most noncommutative'' case, where $\cent(P)$
is a point, corresponds to symplectic manifolds. We will discuss both
situations in the next section.

\begin{remark}\label{rmk:variation} (Variation lattices)

The {\bf regular part} $P_r$ of a Poisson manifold $P$ is the open
dense subset consisting of the regular leaves, i.e. those where the
transverse structure is a zero Poisson structure.  Any Morita
self-equivalence
of $P$ restricts to one on $P_r$, and $\cent(P_r)$ is open
and dense in $\cent(P)$.   Thus, for any $\X \in \Pic(P)$, $h(\X)$
is determined by the restriction to $\cent(P_r)$.
Suppose, now, that the 
symplectic leaf foliation of $P_r$ is a fibration, so that
$\cent(P_r)$ is a manifold.  The proof of Theorem 5.3 of \cite{Xu91}
shows that the restriction of $h(\X)$ to $\cent(P_r)$ preserves the
{\bf variation lattice} of Dazord \cite{Da}, which is a collection of
closed 1-forms attached to integer second homology classes in the fibres,
measuring the variation of the integrals of the fibre symplectic
structures over these classes.  This lattice of closed forms sometimes
provides $\cent(P_r)$ with an additional structure which limits
how $h(\Pic(P))$ may act on $\cent(P)$.  For example, see Section
\ref{sec:liepoisson}.  (Rigidity of $\cent(P_r)$ may also come
from the modular vector field, as in Section \ref{sec:TSS1}.)
\end{remark}

\section{Examples of Picard groups}

\subsection{Symplectic Poisson structures}
We start with a ``trivial'' example.  Any discrete group $G$ is a
symplectic groupoid over a point.  Its symplectic Morita bibundles are
discrete as well, so $\Pic(G)$ for this symplectic groupoid is
isomorphic to $\Outaut(G)$, just as in Example \ref{ex:liegrp}.

Next, suppose that $P=S$ is a connected symplectic manifold.
In this case, $\cent(S)$ is just a point and the exact sequence 
\eqref{eq:exactpois2} is trivial.
The source-simply-connected symplectic groupoid $\Gamma(S)$ is the
fundamental groupoid $\Pi(S)$ (with symplectic structure 
pulled back
from $S\times \overline{S}$), which is transitive since $S$
is connected.  We then have from Example \ref{ex:funda}:

\begin{proposition}
\label{prop:picsymp}
Let $S$ be a connected symplectic manifold. Then
$\Gamma(S)$ is Morita equivalent as a symplectic groupoid to
$\pi_1(S,x)$ for each $x\in S$, and hence
$$
\Pic(S) \cong \Outaut(\pi_1(S,x)).
$$
\end{proposition}
\begin{proof}
It suffices to notice 
 that the Morita bibundle coming from \eqref{eq:gaugebibundle}, 

\begin{equation}\label{eq:bibun}
 \xymatrix {
{\Gamma(S)} \ar@<-0.5ex>[d]\ar@<0.5ex>[d] & \univ{S}\ar[dl]^{\pr}\ar[dr] &
{\pi_1(S,x)} \ar@<-0.5ex>[d]\ar@<0.5ex>[d]\\
S & & \{x\}
 }
\end{equation}
is in fact a {\it symplectic} Morita bibundle, which is clear since the moments are symplectic
realizations.

\end{proof}

As a result of \eqref{eq:inndeg} in Example \ref{ex:gendeg}, we have

\begin{corollary}
If $S$ is a connected symplectic manifold, then $\InnPois(S)$ has the following description:
$$
\InnPois(S) = \{\phi \in \Sympl(S)\; |\; [\phi_*]=\id \in \Outaut(\pi_1(S,x)) \}.
$$
\end{corollary}

So for a connected symplectic manifold, the exact sequence \eqref{eq:exactpoiss1} becomes
\begin{equation}\label{eq:exactsymp}
1 \rightarrow \{\phi \in \Sympl(S)\; |\; [\phi_*]=\id  \in \Outaut(\pi_1(S,x))\}
 \rightarrow \Sympl(S) \stackrel{\inc}{\rightarrow} \Pic(S).
\end{equation}

\begin{remark}
The isomorphism $\eta: \Pic(S) \rightarrow \Outaut(\pi_1(S,x))$
obtained in Proposition \ref{prop:picsymp} can also be seen as a consequence of the
classification of complete Poisson maps with symplectic targets by their holonomy
\cite{SilWein99}.
Recall \cite{SilWein99} that any complete symplectic realization  $M \to S$ is a fibration
with a natural flat connection. A typical fiber, $F$, carries a natural symplectic structure,
preserved under the holonomy action of $\pi_1(S)$. For a choice of base point $x \in S$,
the realization $M \to S$ is isomorphic to $(\univ{S}\times F)/\pi_1(S) \stackrel{\pr}{\to} S$, 
and hence it is completely determined, up to isomorphism, by the holonomy action.

As a consequence, any Morita bimodule $S \stackrel{J_1}{\leftarrow} M \stackrel{J_2}{\rightarrow} S$
is isomorphic to one of the form
$$ 
\begin{diagram}
\node{}\node{\frac{\univ{S}\times F}{\pi_1(S)}}\arrow{sw,t}{\pr}\arrow{se,t}{q_2}\node{}\\
\node{S}\node{}\node{S}
\end{diagram}
$$
The map $q_2$, restricted to $F$, is a covering of $S$, and, since $F$ is simply connected,
there exists a symplectomorphism
$\univ{q_2}:F \rightarrow \univ{S}$ making the following diagram commute:
$$ 
\begin{diagram}
\node{F}\arrow[2]{e,t}{\univ{q_2}}\arrow{se,r}{q_2}\node{}\node{\univ{S}}\arrow{sw,r}{\pr}\\
\node{}\node{S}\node{}
\end{diagram}
$$
The action of $\pi_1(S)$ on $F$ and the map $\univ{q_2}$ induce an action of $\pi_1(S)$
on $\univ{S}$ by deck transformations, and this action defines an automorphism $\phi$
of $\pi_1(S)$. One can check that 
$$
\eta([M]) = [\phi] \in \Outaut(\pi_1(S)).
$$
\end{remark}

\begin{remark}
Gompf \cite{Gompf95} has shown that every finitely presented group is
the fundamental group of a compact symplectic 4-manifold.  (Without the compactness
and dimension restriction, it is much easier to realize these groups
as fundamental groups of cotangent bundles.)  This shows that the
subcategories of $\catpoiss$ consisting of groups and of symplectic
manifolds are essentially the same as far as their objects are
concerned.  Comparing them on the level of morphisms means deciding
which outer automorphisms of the fundamental group are realizable by
symplectic diffeomorphisms.  (For closed surfaces, Nielsen's theorem
on mapping classes \cite{Nielsen} and Moser's theorem on volume elements \cite{Moser} 
together imply that all outer automorphisms are realizable by symplectic or
antisymplectic diffeomorphisms.)
\end{remark}

\subsection{The zero Poisson structure}

Let $(P,\pi)$ be a Poisson manifold with $\pi=0$. In this case,
 $\cent(P) = P$,
 $\Aut(P) = \Diff(P)= \Poiss(P)$, and the inner Poisson automorphisms
are trivial. The exact sequence \eqref{eq:exactpoiss1} becomes
\begin{equation}\label{eq:exinc}
1 \to \Poiss(P) \stackrel{\inc}{\rightarrow} \Pic(P).
\end{equation}
For \eqref{eq:exactpois2}, we get 
\begin{equation}\label{eq:exproj}
1 \to \PicZP(P) \to \Pic(P) \stackrel{\proj}{\rightarrow} \Poiss(P),
\end{equation}
and, as in the case of commutative algebras, $\proj \circ \inc = \id$.

Let $\pr:T^*P \to P$ be the natural projection, and let $\omega$ 
denote the canonical symplectic form on $T^*P$. The following lemma
is closely related to \cite[Thm.~4.1]{Naga} and follows from \cite[Prop.~4.7.1]{Wood92}.

\begin{lemma}\label{lem:cot}
Let $(P,\pi)$ be a Poisson manifold with $\pi =0$.
Let $\omega'$ be a symplectic form on $T^*P$
for which the
$\pr$-fibres and the zero section
$P \hookrightarrow T^*P$ are lagrangian submanifolds and
$\pr$ is a complete Poisson map. Then there exists a symplectomorphism
$f:(T^*P,\omega) \to (T^*P,\omega')$ so that $f(P)=P$ and $f\circ \pr = \pr \circ f$.
\end{lemma}

We now give an explicit description of $\PicZP(P)$ and $\Pic(P)$.

\begin{proposition}\label{prop:piczero}
Let $(P,\pi)$ be a Poisson manifold with $\pi =0$. Then
\begin{itemize}
\item[(i)] $\PicZP(P) \cong H^2(P,\mathbb{R})$;
\item[(ii)] $\Pic(P) \cong  \Diff(P) \ltimes H^2(P,\mathbb{R})$, where the semi-direct product
is with respect to the natural action of diffeomorphisms on cohomology by pull-back.
\end{itemize}
\end{proposition}

\begin{proof}
The canonical source-simply-connected symplectic groupoid of $P$ is 
$T^*P$, equipped with its canonical symplectic form $\omega$
and groupoid structure given by fibrewise addition.

It is easy to check that, for any closed 2-form
$\alpha \in \Omega^2(P)$, $(T^*P, \omega + \pr^*\alpha)$ is a 
$(P,P)$-Morita bimodule defining an element in $\PicZP(P)$.
Let us consider the map
\begin{equation}
\tau: \{\alpha \in \Omega^2(P), \; \alpha \mbox{ closed}\,\} \to \PicZP(P), \;\;\;  
\alpha \mapsto [(T^*P, \omega + \pr^*\alpha)].
\end{equation}

\noindent{\it Claim 1: The map $\tau$ is onto.} 

In order to prove Claim 1, note that, by Example~\ref{ex:trivcontr} (see \eqref{eq:picztriv}),
any $(P,P)$-Morita bimodule $S$ inducing the identity map on $P$ must be of the form
$(T^*P,\nu)$, where $\nu$ is a symplectic form for which the $\pr$-fibres are lagrangian and
$\pr$ is complete. Consider the closed 2-form $\beta=\iota^*\nu \in \Omega^2(P)$, where $\iota:P \to T^*P$
is the zero-section embedding. Then $\omega' = \omega -\pr^*\beta $ satisfies the conditions
of Lemma \ref{lem:cot}, so there is a symplectic bimodule isomorphism 
$f: (T^*P,\omega) \to (T^*P,\omega').$ But $f^*\omega' = \omega$ implies that
$$
f^*\omega - f^*\pr^*\beta = f^*\omega - \pr^*f^*\beta = \omega.
$$
So $f^*\omega = \omega + \pr^*(f^*\beta)$, and therefore $[S]= \tau(\alpha)$, where
$\alpha = f^*\beta$.

\noindent{\it Claim 2: $\Ker(\tau) = \{\alpha \in \Omega^2(P), \; \alpha \mbox{ exact}\,\}$.}

If $\alpha \in \Ker(\tau)$, then there exists a diffeomorphism $f:T^*P \to T^*P$ with
$f^*(\omega + \pr^*\alpha)= \omega$. Since $\omega$ is exact, so is $\alpha$.
On the other hand, if $\alpha = d\theta$, for $\theta \in \Omega^1(P)$, then 
the map $f:T^*P \to T^*P$, $(x,\xi) \mapsto (x, \xi - \theta)$ satisfies
$f^*(\omega + \pr^*d\theta) = \omega$ and defines a symplectic bimodule isomorphism.
So $\alpha \in \Ker(\tau)$.

It is not hard to check that $\tau$ is a group homomorphism. So
we have an exact sequence of groups
$$
1 \to \{\alpha \in \Omega^2(P), \; \alpha \mbox{ exact}\,\} \to
\{\alpha \in \Omega^2(P), \; \alpha \mbox{ closed}\,\} \stackrel{\tau}{\to} \PicZP(P) \to 1,
$$
so $\PicZP(P) \cong H^2(P,\mathbb{R})$. This proves $(i)$.

Since the maps $\inc$ and $\proj$ in \eqref{eq:exinc} and \eqref{eq:exproj} 
satisfy $\proj \circ \inc=\id$, we have an identification
$$
\Diff(P) \times H^2(P,\mathbb{R}) \to \Pic(P), \;\; 
(\phi, \alpha) \mapsto ([(T^*P)_{\phi}], \omega + \pr^*\alpha),
$$
under which the tensor product on $\Pic(P)$ becomes
$$
((\phi, \alpha) , (\psi, \beta)) \mapsto (\phi\circ \psi, \alpha + (\phi^{-1})^*\beta).
$$
This proves $(ii)$.
\end{proof}

\begin{remark}
As Dmitry Roytenberg has suggested to us, one can pass directly from
an element in $\PicZP(P)$ to a 1-dimensional extension of the tangent
bundle Lie algebroid $TP$.  Such extensions are classified by
elements of $H^2(P,\reals)$; a vector bundle splitting of the extension has a
``curvature'' whose cohomology class determines the extension up to
isomorphism.
  
Let $S$ be a $(P,P)$-Morita bimodule with maps $J_1=J_2=J$.
The fibres of $J$ are the leaves of a
lagrangian foliation, hence they carry flat affine structures.  
For each $x \in P$, the affine
functions on $J^{-1}(x)$ form a vector space which is a fibre of the
vector bundle $E$ over $P$ whose sections are the fibrewise affine
functions on $S$.  Given any such section, its fibre derivative may be
identified with a vector field on $P$; the resulting bundle map from
$E$ to $TP$ is the anchor map for a Lie algebroid structure in which
the bracket is the Poisson bracket of $S$, restricted to fibrewise
affine functions.  Splittings of this Lie algebroid extension
correspond to cross sections of $J$; this correspondence may be used
to show that the characteristic class of the Lie algebroid is the same
as that attached to the bimodule in Proposition \ref{prop:piczero}.
\end{remark}

The Picard group of $P$ is the same as that of 
the source-simply-connected symplectic groupoid $\Gamma(P)=T^*P$.  
However, there are other symplectic
groupoids over $P$, for which the calculation of the Picard group is
more complicated.  

We have already seen that, when $P$ is a point,
dropping the source-connectivity assumption leads to the class of
discrete groups as symplectic groupoids.  In general, the symplectic 
groupoids over $P$ for which the (components of the)
source fibres have trivial
fundamental groups consist of action groupoids for the actions of discrete groups on
$T^*P$ which preserve the foliation by the fibres.  These actions
may be described as ``affinizations'' of the standard cotangent lifts of
actions on $P$.  

On the other hand, we may keep source-connectivity but allow
source fibres to have nontrivial fundamental groups.  
In this case, the groupoid is obtained by
dividing the cotangent bundle $T^*P$ by a ``lattice'' $\Lambda$ which
is a bundle of discrete subgroups with local bases consisting of
 $k$-tuples of 
closed $1$-forms for some $k \leq \dim (P)$.  Such lattices arise
naturally in the theory of action-angle variables as analyzed by
Duistermaat \cite{Duis}.  The lattice describes what Duistermaat
refers to as the monodromy of a lagrangian fibration.  Once the
monodromy (i.e. the groupoid) has been fixed, then the integrable
systems with this monodromy are precisely the symplectic bibundles for
this groupoid; their isomorphism classes are the Picard group, with
respect to the ``tensor product'' operation.   

Analysis of this Picard
group begins as follows \cite{Duis,Zung}.  When the monodromy is
''trivial,'' i.e. $\Lambda$ has a global basis of sections, the next
invariant of a torus fibration is its Chern class, which is an element
of $H^2(P)$ with values in the sections of $\Lambda.$ 
Finally, when the Chern class is trivial, there
is a {\bf lagrangian class} in $H^2(P,\reals)$ whose vanishing is the condition
for the fibration to admit a {\em lagrangian} cross section.  In
general, these invariants are mixed in a complicated way (for
instance, the Chern class lies in the cohomology of a sheaf determined
by the monodromy). The starting point for the description of the Picard
group in terms of these invariants is Example \ref{ex:abelian} combined
with an extension of Lemma \ref{lem:cot} to more general 
lagrangian foliations. We will leave this discussion to a separate paper.

\subsection{Lie--Poisson structures on duals of compact groups}
\label{sec:liepoisson}
For any Lie group $G$, the
cotangent bundle $T^*G$ is a symplectic
groupoid over the dual Lie algebra $\gstar$; it is
source-simply-connected just when $G$ is (connected and)
simply-connected.  What is the Picard group of $T^*G$?

When $G$ is a torus, we
are in the situation described after Propostion~\ref{prop:piczero}.  
In the remainder of this section, we will discuss the case where $G$
is compact and simply-connected, so that we are also studying the Picard
group of $\gstar$.  Here, the regular part $\gstar_r$ is
an open dense subset, in which the leaves are all simply connected, 
and the leaf space
$\cent(\gstar_r)$ may be identified with the interior of a Weyl
chamber $W$.  $\gstar_r$ itself is a (trivial) bundle of symplectic manifolds
over $\cent(\gstar_r)$ whose fibres are flag manifolds $G/T$, with
symplectic structure depending on the base point in $W$.  By Remark
\ref{rmk:variation} above, for any $\X\in \Pic(\gstar)$
the restriction to $\cent(\gstar_r)$ of its action $h(\X)$ on the leaf
space must preserve the variation lattice.
On the other hand, this variation lattice is well known to be the
weight lattice, once one identifies the tangent spaces to $W$ with the
vector space $\frakt^*$.  Since the forms making up the variation
lattice are
constant with respect to the natural affine structure of the Weyl
chamber, $h(\X)$ must be (the restriction of) a linear map.  Linear
maps of this type, preserving the Weyl chamber and the weight lattice,
arise from permutations of the fundamental weights.  It follows from 
Xu's work \cite{Xu91} 
that all these permutations arise from the Picard group of
$\gstar_r$.  On the other hand, for the full $\gstar$, one may
restrict a Picard group element to the fibre of $T^*G$ over the zero
element of $\gstar$, obtaining an outer automorphism of $G$.  These
correspond to automorphisms of the Dynkin diagram, which are generally
quite special compared to the permutations of the fundamental weights.
Taking into account the simple-connectivity of the symplectic leaves,
we believe that this homomorphism from $\Pic(\gstar)$ to $\Outaut(G)$
may in fact be an isomorphism.  

\section{Gauge versus Morita equivalence}

We recall here the  notion of a gauge transformation of a Poisson
structure associated with a 
closed $2$-form \cite{SeWe01}.

Let $(P,\pi)$ be a Poisson manifold, and let $B \in \Omega^2(P)$ be a closed
$2$-form. We identify $\pi$ (resp.~$B$) with the bundle map
$\widetilde{\pi}:T^*P \rightarrow TP$, $\widetilde{\pi}(\alpha)(\beta) = \pi(\beta,\alpha)$
(resp.~$\widetilde{B}:TP \rightarrow T^*P$, $\widetilde{B}(u)(v) = B(u,v)$).
If the bundle endomorphism
\begin{equation}\label{eq:bundlemap}
1 + \widetilde{B} \widetilde{\pi}: T^*P \rightarrow T^*P
\end{equation}
is invertible, we define the {\bf gauge transformation} $\tau_B$ of $\pi$ associated with
$B$ by
\begin{equation}\label{eq:transform}
\widetilde{\tau_B(\pi)}=  \widetilde{\pi}(1 + \widetilde{B}\widetilde{\pi})^{-1}.
\end{equation}
When $\pi$ is nondegenerate, \eqref{eq:transform} just says that 
$\widetilde{\tau_B(\pi)}^{-1}=\widetilde{\pi}^{-1} + B$; 
in particular, any two symplectic structures
on a given manifold are gauge equivalent.
In general, a gauge transformation $\tau_B$ produces a Poisson
 structure
 whose leaf decomposition
is the same as before, but the symplectic structures along the leaves
 differ by the pull-backs 
of $B$.  
\footnote{Since $1+\widetilde{B}\widetilde{\pi}$ might not be
 invertible, we are not quite dealing here with a group action, but
 with a groupoid, obtained by restricting to the Poisson structures
 the groupoid associated with the action of the additive group of
 closed $2$-forms on the space of {\bf Dirac structures} \cite{Cou90}
 on $P$ (see \cite{SeWe01}).}
We call two Poisson structures $\pi_1$ and $\pi_2$ on $P$ {\bf gauge
equivalent} if there exists a closed $2$-form $B$ with
$\tau_B(\pi_1)=\pi_2$.  More generally, two Poisson manifolds
$(P_1,\pi_1)$ and $(P_2,\pi_2)$ are called {\bf gauge equivalent up to
Poisson diffeomorphism} if there exists a Poisson diffeomorphism
$f:(P_1,\pi_1) \rightarrow (P_2,\tau_B(\pi_2))$ for some closed
$2$-form $B \in \Omega^2(P_2)$.

The close relationship between gauge and Morita equivalences is 
shown by the following result \cite{BuRad02}.

\begin{theorem}\label{thm:gaugemorita}
If two integrable Poisson manifolds $(P_1,\pi_1)$ and $(P_2,\pi_2)$ are gauge equivalent
up to Poisson diffeomorphism, then they are Morita equivalent. The converse does not hold in general.
\end{theorem}

The proof of the theorem follows from the description of the effect of gauge transformations on symplectic
groupoids. Let $(P,\pi)$ be an integrable  Poisson manifold, and
let $\Gamma(P)$ be its source-simply-connected symplectic groupoid.
Since $\pi$ and any gauge equivalent 
Poisson structure $\tau_B(\pi)$ have isomorphic Lie algebroids,
$\tau_B(\pi)$ is integrable, and its (source-simply-connected)
symplectic groupoid can be identified, as a Lie groupoid, 
with $\Gamma(P)$. 
The source and target maps are unchanged, but the original symplectic
structure $\omega$ on $\Gamma(P)$ is changed to 
$\tau_B(\omega):= \omega + t^*B -s^*B.$  (The additional terms form a
coboundary in the ``bar-de Rham'' double complex attached to a
groupoid, as discussed, for instance, in \cite{BeXuZh}.  
The meaning of this fact is not very clear to us, but see \cite[Ex.~6.6]{BCWZ}
and Section \ref{sec:piclie}.)

In other words, we have the correspondences: 
\begin{equation}\label{eq:gaugegroupoid}
 \xymatrix {
{(\Gamma(P),\omega)} \ar@<-1ex>[d]_{t}\ar@<1ex>[d]^{s} \ar@{.>}[r] & 
{(\Gamma(P), \tau_B(\omega))} \ar@<-1ex>[d]_{t}\ar@<1ex>[d]^{s}\\
(P,\pi) \ar@{.>}[r]^-{\tau_B} & (P,\tau_B(\pi))
 }
\end{equation}
A Morita bimodule for $\pi$ and $\tau_B(\pi)$ is obtained by means of a ``half'' twist:
\begin{equation}
 \xymatrix@!0@=1.5cm {
& {(\Gamma(P),\omega')} \ar@<0ex>[dr]^{s} \ar@<0ex>[dl]_{t} &\\
(P,\pi) & & (P,\tau_B(\pi))
 }
\end{equation}
where $\omega'=\omega-s^*B$.

To see that Morita equivalence does not imply gauge equivalence, even
up to Poisson diffeomorphism, we
note that two symplectic manifolds are gauge equivalent
up to Poisson diffeomorphism if and only if they are symplectomorphic, whereas
Morita equivalence just amounts to isomorphism of fundamental
groups. Example 5.2 in \cite{BuRad02} shows that even Morita
equivalent Poisson structures on the {\it same} manifold can fail to 
be gauge equivalent
up to Poisson diffeomorphism.  On the other hand, we will present
in Section \ref{sec:TSS1} a class of Poisson structures on surfaces
for which the notions of gauge and Morita equivalence do coincide.

\begin{remark}
As we mentioned above, the gauge transformation construction extends
to an action of the closed $2$-forms on $P$ on the set of Dirac
structures on $P$.  If $(P,\pi)$ is an integrable Poisson manifold
with symplectic groupoid $(\Gamma(P),\omega)$ such that $\tau_B(\pi)$
fails to be a Poisson structure, then the $2$-form $\omega + t^*B
-s^*B$ becomes degenerate, and the diagram \eqref{eq:gaugegroupoid}
suggests that Dirac structures should be ``integrated'' to groupoids
carrying multiplicative pre-symplectic forms, generalizing the
correspondence between integrable Poisson manifolds and symplectic
groupoids. The precise relationship between Dirac structures and
``pre-symplectic groupoids'' is developed in \cite{BCWZ}, where the
situation of Dirac structures ``twisted'' by closed $3$-forms
\cite{SeWe01} is also investigated.  Following \cite{Xu}, one may use the
correspondence
$$
\mbox{ (twisted) Dirac structures } \leftrightarrow
\mbox{ (twisted) pre-symplectic groupoids}
$$
to give a groupoid interpretation of
quasi-hamiltonian
actions and their group-valued momentum maps \cite{AMM98}.
\end{remark}

\begin{remark}
The results in \cite{RiefSch}, on Morita equivalence of quantum tori,
and \cite{BuWa2001,JSW2001}, on Morita equivalence of formal
deformation quantizations,
suggest that one should single out for special attention the gauge
transformations by closed $2$-forms which belong to integer cohomology
classes.  It should be interesting to investigate how 
the geometric properties of
the symplectic bimodules given by Theorem \ref{thm:gaugemorita} in
this case might relate to the Morita equivalence on the quantum level.
\end{remark}

\section
{Topologically stable structures on surfaces}\label{sec:TSS1} 

This section has been written in collaboration with Olga Radko.
 
Let \( \Sigma  \) be a compact connected oriented surface equipped 
with a Poisson structure \( \pi  \) which has at most linear 
degeneracies.  We call 
such a structure {\bf topologically stable}, since the topology of 
its zero set is preserved under small perturbations of $\pi$.  
This zero set consists of 
$n$ smooth disjoint, closed curves on \( \Sigma  \), for some $n\geq 
0$.  Each of them carries a natural orientation given by any 
modular vector field for $\pi$.  We denote by $Z(\Sigma ,\pi )$ 
the zero set, considered as an oriented 1-manifold.

We call two topologically stable surfaces \( (\Sigma ,\pi ) \) and \(
(\Sigma ',\pi ') \) {\bf topologically equivalent} if there is an
orientation-preserving diffeomorphism \( \varphi :\Sigma \to \Sigma '
\) such that \( \varphi (Z(\Sigma ,\pi ))=Z(\Sigma ',\pi ') \).  The
associated equivalence class is denoted by \( [Z(\Sigma ,\pi )] \).
 
Each topologically stable structure \( (\Sigma ,\pi ) \) with 
\( n \) zero curves has \( (n+1) \) numerical invariants: 
\( n \) modular periods 
(periods of a modular vector field around the zero curves) and a regularized volume invariant 
(generalizing the Liouville volume in the symplectic case). Together 
with the topological equivalence class, these invariants completely 
classify topologically stable structures up to orientation-preserving 
Poisson isomorphisms 
\cite{Rad01}. 
 
A topological equivalence class $[Z(\Sigma,\pi)]$ can be encoded by an
oriented labeled graph $\mathfrak{G}(\Sigma,\pi)$, with
a vertex for  
each \( 2 \)-dimensional leaf of the structure, two vertices being
connected by an edge for each boundary  
zero curve which they share.  Each edge is oriented so that it points 
toward the vertex for which the Poisson structure is positive with 
respect to the orientation of $\Sigma$.   
We label each vertex by the genus 
of the corresponding leaf.  (Note that this genus, together with the 
number of edges at the vertex, completely determines the topology of 
the leaf.)  If we also assign to each edge of \( \mathfrak
{G}(\Sigma ,\pi ) \)  
the modular period of \( (\Sigma ,\pi ) \) around the corresponding 
zero curve, we obtain a more elaborately labeled graph which we denote
	      by \( \mathfrak {G}_{T}(\Sigma ,\pi ) \).

\subsection{Morita equivalence of topologically stable structures}

Topologically stable structures (TSS) on surfaces form a class 
of Poisson structures for which Morita equivalence 
is the same as the gauge equivalence up to a diffeomorphism. The first
step toward establishing this is:

\begin{theorem}
Two TSS \( (\Sigma _{1},\pi _{1}) \) and \( (\Sigma _{2},\pi _{2}) \)
are Morita equivalent iff there is an isomorphism of labeled graphs
\( \mathfrak {G}_{T}(\Sigma _{1},\pi _{1})\simeq \mathfrak {G}_{T}(\Sigma _{2},\pi _{2}) \). 
\end{theorem}
\begin{proof} 

Assume that \( (\Sigma _{1},\pi _{1}) \) and \( (\Sigma _{2},\pi _{2})
\) are Morita equivalent.  As was shown in \cite{BuRad02} for the case
of the sphere, the fact that Morita equivalence implies homeomorphism
of the leaf spaces means that \( \mathfrak {G}(\Sigma _{1},\pi
_{1})\simeq \mathfrak {G}(\Sigma _{2},\pi _{2}) \).  Since the
restriction of the Morita equivalence bimodule to each leaf is a
Morita equivalence (see Remark \ref{remark:restrict}, it follows that
the genera of the associated leaves are the same. Finally, since
modular periods are invariants of Morita equivalence of topologically
stable structures \cite{BuRad02}, \( \mathfrak {G}_{T}(\Sigma _{1},\pi
_{1})\simeq \mathfrak {G}_{T}(\Sigma _{2},\pi _{2}) \) as labeled
graphs.

Conversely, if \( \mathfrak {G}_{T}(\Sigma _{1},\pi _{1})\simeq
  \mathfrak {G}_{T}(\Sigma _{2},\pi _{2}) \), 
there is a diffeomorphism $\alpha: \Sigma _{1}\to \Sigma _{2}$ which
 maps $Z(\Sigma_1,\pi_1)$ to $Z(\Sigma_2,\pi_2)$, 
since the graph \( \mathfrak
{G}(\Sigma ,\pi ) \) completely encodes the topology of the decomposition of \(
\Sigma \) into its \( 2 \)-dimensional symplectic leaves. On
$\Sigma=\Sigma _{1}$,
let \( \pi =\pi _{1},\, \pi '=\alpha _{*}^{-1}(\pi _{2}) \).
Since \( \mathfrak {G}_{T}(\Sigma _{1},\pi _{1})\simeq \mathfrak
{G}_{T}(\Sigma _{2},\pi _{2}) \), we know that, for any zero curve \(
\gamma \in Z(\pi )=Z(\pi ') \), the modular periods for the two
structures are the same.  By Theorem 6.2 of \cite{BuRad02}, $\pi$ and
  $\pi'$ are gauge equivalent up to Poisson diffeomorphism, and so are 
 $(\Sigma _{1},\pi _{1}) \) and \( (\Sigma _{2},\pi _{2}) $; 
hence they are Morita equivalent by Theorem \ref{thm:gaugemorita}.
\end{proof}

The fact that gauge equivalence up to a diffeomorphism implies Morita
equivalence and the second part of the proof above imply:

\begin{corollary}
For TSS the notions of Morita and gauge equivalence up to a diffeomorphism
coincide.\end{corollary}

\subsection{Picard groups of topologically stable structures}\label{sec:TSS2} 

In this section, we merely raise the problem of computing the
Picard group of a TSS $(\Sigma,\pi)$.  It appears from the arguments
above that the ingredients of $\Pic(\Sigma,\pi)$ should be:

\begin{enumerate}
\item the automorphism group of the labeled graph
\( \mathfrak {G}_{T}(\Sigma ,\pi ) \);

\item the torus which is the product of the groups of rotations of the
zero curves;

\item the  outer automorphism groups of the fundamental groups of
  the $2$-dimensional leaves.

\end{enumerate}

How these ingredients are combined should be described by an
algebraic/combinatorial object which is a further refinement of \(
\mathfrak {G}_{T}(\Sigma ,\pi ) \), and which encodes the inclusions
of the fundamental groups of the zero curves into those of the
adjacent $2$-dimensional leaves. 

\section{Further questions}
We hope that this paper represents the beginning of an interesting
line of research.  With that in mind, we conclude with discussion of
some issues which remain at least partly unresolved.

\subsection{The Lie algebra of the Picard group}
\label{sec:piclie}
In all the examples above, the Picard group consists of a
discrete part associated with fundamental groups, and a continuous
part associated with diffeomorphisms of leaf spaces.  It seems useful,
therefore, to think of Picard groups as  (possibly
infinite-dimensional) Lie groups, and to study the infinitesimal
objects which should play the role of their Lie algebras.

It should be possible to define {\bf Picard Lie algebras}
$\Piclie(\A)$, $\Piclie(\Gamma)$, and $\Piclie(P)$ for algebras, (Lie,
symplectic) groupoids, and Poisson manifolds.  Their elements should
be isomorphism classes of infinitesimal deformations of the identity
bimodule (or bibundle), with the underlying object held fixed.  (This
is to be contrasted with the analysis in \cite{BuWa} of deformations of
bimodules as the underlying algebra is deformed.)  In the case of
algebras, if one fixes the underlying $k$-module of the bimodule to be
$\A$ itself, $\Piclie(\A)$ has a description in terms of
Hochschild cocycles and Gerstenhaber brackets.  For symplectic
groupoids, one may fix the underlying manifold of the bibundle and
deform the symplectic structure along with the bibundle structure.
The former deformations should be related to the groupoid double
complex mentioned after Theorem \ref{thm:gaugemorita} and produce 
2-cohomology classes on the base, while the 
latter should lead to the ``diffeomorphism'' part of $\Piclie(\Gamma)$.

Revisiting the exact sequences in Section \ref{sec:picsymp} and
passing to Lie algebras, we obtain first of all from
\eqref{eq:exactpoiss1} and \eqref{eq:exactbis} the sequence
$$
0\to C^{\infty}(P)/H_\pi^0(P)\stackrel{d_\pi}{\to} \Chi_{\pi}(P) \to \Piclie(P),
$$ 
so that $H_\pi^1(P)$ sits inside $\Piclie(P)$ as a subalgebra.
(Here $\Chi_{\pi}(P)$ denotes the space of Poisson vector fields,
$H^{\bullet}_{\pi}(P)$ are Poisson cohomology spaces and $d_{\pi}$ is
the Poisson differential, see e.g. \cite{SilWein99}.)
Similarly,
\eqref{eq:exactpois2} becomes
$$
0\to \Piclie_{\mathcal Z (P)} \to \Piclie(P)\to \Chi(\mathcal Z
(P)).
$$
Here, $\Chi(\mathcal Z (P))$ denotes the space of vector fields on
$\mathcal Z (P), $ restricted when appropriate to those whose flows
preserve an additional structure such as a variation lattice.

\subsection{Representation equivalence vs. Morita equivalence}

Although, for algebras, representation equivalence implies Morita
equivalence, this is not the case in other categories, where an object
may not serve as its own identity bimodule.  Landsman \cite{Land00b}
discusses this phenomenon, including
the case of $C^*$-algebras, 
where the  notion of representation must be suitably chosen
so that representation equivalence does imply Morita equivalence 
\cite{Ri74}.   

For Poisson manifolds, Xu gave an explicit counterexample in \cite{Xu91}.  
One simply takes a Poisson manifold and multiplies the Poisson structure by $2$.
This clearly does not change the category of representations, but
examples may be given in which the Morita equivalence class does
change, since the variation lattice on the leaf space is multiplied by
$2$.  (See Remark \ref{rmk:variation}.)
In this setting, it seems that, rather than changing the notion of
representation, we should give a refined structure to the {\em category}
of representations.  For instance, in Xu's examples, one may distinguish
the two representation categories by considering the Lie algebras of
the automorphism groups as (codimension 1 quotients of) Poisson algebras. 
An interesting approach would be to define a notion of representation
for Poisson manifolds which is refined enough to detect the
difference between the two examples, but we will reserve this discussion
for a future work.

\subsection{Nonintegrable Poisson manifolds}\label{sec:extended}

If we try to build a category in which the Morita equivalences are
isomorphisms, and which contains {\em all} Poisson manifolds and complete
Poisson maps, two problems immediately arise.  First of all, if we
define morphisms to be (symplectic) bimodules, then only the integrable
Poisson manifolds admit identity morphisms.  Second, although we can
associate a bimodule to any complete Poisson map between integrable
Poisson manifolds, these bimodules are generally not generalized
morphisms (see Remark \ref{rem:etale}), so that
it becomes difficult to define compositions.  Both of these problems
occur because certain natural constructions lead to
leaf spaces of foliations, which may not be manifolds.

It seems that the difficulties just described 
can be circumvented if one admits as symplectic bimodules not only
symplectic manifolds, but also the leaf spaces of transversely
symplectic foliations.  
The technical aspects of this approach are yet to be worked out,
requiring an appropriate category of ``leaf spaces,'' so we merely
sketch some ideas here.

The problem of identity morphisms is handled as follows.  Cattaneo and
Felder \cite{CaFe} constructed, for every Poisson manifold $P$, a groupoid
$\Gamma(P)$ over $P$ which they identify as the ``phase space'' for
a Poisson sigma model associated with $P$.  More concretely,
$\Gamma(P)$ is the quotient of a Banach manifold of paths in $T^*P$
by a foliation of finite codimension (twice the dimension of $P$).
Crainic and Fernandes \cite{CrFe01} extended this construction to the
case where the Poisson manifold $P$ are
replaced by an arbitrary manifold $M$ and a Lie algebroid $A$ over
$M$ (the Lie algebroid in the Poisson case being the cotangent bundle
$T^*P$), obtaining a canonical groupoid $\Gamma(A)$ over $M$ which
they call the ``Weinstein groupoid'' of $A$.  
The ideas of this construction have appeared independently in
\cite{Severa}. 
By abuse of notation, we
will denote $\Gamma(T^*P)$ for a Poisson manifold $P$ by
$\Gamma(P)$.

Since the elements of $\Gamma(A)$ are homotopy classes of paths
(in a sense adapted to the presence of $A$) in $M$, we propose, 
following \cite{Severa}, to call $\Gamma(A)$
the {\bf fundamental groupoid} of $(M,A)$ or, when $A$ is the
cotangent Lie algebroid $T^*P$ of a Poisson manifold $P$, the {\bf
Poisson fundamental groupoid} of $P$.  Notice that $\Gamma(A)$ is
the usual fundamental groupoid of $M$ just when $A$ is isomorphic to
$TM$ (i.e. when $P$ is symplectic in the Poisson case).

Although $\Gamma(A)$ it is generally not a manifold, it is always
the leaf space of a foliation.  This foliation originally lives in the
Banach manifold of $A$-paths (see \cite{CrFe01}), but it is
finite-codimensional, so if 
one prefers to work with finite dimensional spaces, one may identify 
the leaf space with the
the orbit space of an \'etale groupoid obtained by restricting the
holonomy groupoid  of the foliation to a complete transversal.  Crainic
and Fernandes gave explicit criteria in \cite{CrFe01}
for $\Gamma(A)$ to be an ordinary manifold, with 
refinements in \cite{CrFe02} for the Poisson case.

Leaf spaces may still be treated for many purposes as if it they were
manifolds.  The best way to do this seems to define a
{\bf leafspace structure} on 
$X$ as a Morita equivalence class of \'etale groupoids
whose orbit spaces are identified with $X$.  Mappings between
leafspaces are then equivalence classes of regular bimodules over these
groupoids, and a symplectic structure on $X$ is 
given by  groupoid-invariant
symplectic structures on the manifolds of objects of  \'etale groupoids
realizing $X$.\footnote{Note that the term ``symplectic
  leaf spaces'' has quite a different meaning.}
All this analysis is perhaps best carried out in the language of
``differentiable stacks'' \cite{Pr}.

We may now define a category in which the objects are Poisson
manifolds, and in which a morphism 
 $P_2\leftarrow P_1$ is an
isomorphism class of left principal bibundles
of the form \eqref{eq:symbim} in which $S$ is a symplectic leafspace.
An isomorphism in this category will still be called a Morita
equivalence.  Fortunately,
it follows from the local triviality of principal bundles \cite{MM03}
that the total space of a Morita
equivalence between {\em integrable} Poisson manifolds must be a
smooth manifold, so that the manifolds are already Morita equivalent
in the usual sense. 

\begin{footnotesize}

\end{footnotesize}


\begin{thebibliography}{10}

\bibitem {AlDa91}
{\sc Albert, C.,  Dazord, P.: }\newblock {\em Groupo\"{\i}des de Lie et groupo\"{\i}des
symplectiques}.
\newblock In: {\em Symplectic geometry, groupoids,   
and integrable systems, S\'{e}minaire sud-Rhodanien 
de g\'{e}om\'{e}trie \`{a} Berkeley (1989)}, P. Dazord and A.
Weinstein, eds., 1--11. Springer-MSRI Series, 1991.

\bibitem {AMM98}
{\sc Alekseev, A., Malkin, A., Meinrenken, E.: }\newblock {\em Lie group valued
  moment maps}.
\newblock J. Differential Geom.  {\bf 48} (1998), 445--495.

\bibitem {Bass68}
{\sc Bass, H.: }\newblock {\em Algebraic ${K}$-theory}.
\newblock W. A. Benjamin, Inc., New York-Amsterdam, 1968.

\bibitem{BeXuZh}
{\sc Behrend. K.,  Xu, P.,  Zhang, B.: }
\newblock {\em Equivariant gerbes over compact simple Lie groups}.
\newblock to appear in C. R. Acad.~Sci.~Paris.

\bibitem {Ben67}
{\sc B{\'e}nabou, J.: }\newblock {\em Introduction to bicategories}.
\newblock In: {\em Reports of the Midwest Category Seminar},   1--77. Springer,
  Berlin, 1967.

\bibitem {BCWZ}
{\sc Bursztyn, H., Crainic, M., Weinstein, A., Zhu, C.: }\newblock {\em
  Integration of twisted Dirac brackets}. \newblock Math.DG/0303180.  

\bibitem {BuRad02}
{\sc Bursztyn, H., Radko, O.: }\newblock {\em Gauge equivalence of Dirac
  structures and symplectic groupoids}.
\newblock Ann. Inst. Fourier (Grenoble) {\bf 53} (2003), 309--337.

\bibitem {BuWa2001}
{\sc Bursztyn, H., Waldmann, S.: }\newblock {\em The characteristic classes of
  {M}orita equivalent star products on symplectic manifolds}.
\newblock Comm. Math. Phys.  {\bf 228} (2002), 103--121.

\bibitem {BuWa}
{\sc Bursztyn, H., Waldmann, S.: }\newblock {\em Bimodule deformations, Picard groups
and contravariant connections}. \newblock Math.QA/0207255, to appear
in $K$-Theory.  

\bibitem {SilWein99}
{\sc Cannas~da Silva, A., Weinstein, A.: }\newblock {\em Geometric models for
  noncommutative algebras}.
\newblock American Mathematical Society, Providence, RI, 1999.

\bibitem {Ca2001}
{\sc Cartier, P.: }\newblock {\em  A mad day's work: from Grothendieck to Connes and
Kontsevich. The evolution of concepts of space and symmetry}.
Translated from the French by Roger Cooke.
\newblock Bull. Amer. Math. Soc. (N.S.)  {\bf 38}.4 (2001), 389--408.

\bibitem {CaFe}
{\sc Cattaneo, A., Felder, G.: }\newblock {\em Poisson sigma models and symplectic groupoids}.
\newblock In: {\em Quantization of singular symplectic quotients},   
Progr. Math. 198, 61--93, Birkh\"auser, Basel, 2001.

\bibitem {CDW87}
{\sc Coste, A., Dazord, P., Weinstein, A.: }\newblock {\em Groupo\"\i des
  symplectiques}.
\newblock In: {\em Publications du D\'epartement de Math\'ematiques. Nouvelle
  S\'erie. A, Vol.\ 2},   i--ii, 1--62. Univ. Claude-Bernard, Lyon, 1987.

\bibitem {Cou90}
{\sc Courant, T.: }\newblock {\em Dirac manifolds}.
\newblock Trans. Amer. Math. Soc.  {\bf 319} (1990), 631--661.

\bibitem {Cra} 
{\sc Crainic, M.:}\newblock {\em Differentiable and
algebroid cohomology, van Est isomorphisms, and characteristic
classes}.   \newblock Math.DG/0006064, to appear in Comm. Math. Helv. 

\bibitem {CrFe01}
{\sc Crainic, M., Fernandes, R.: }\newblock {\em Integrability of Lie
  brackets}. \newblock Math.DG/0105033, to appear in Annals of Math. 

\bibitem {CrFe02}
{\sc Crainic, M., Fernandes, R.: }\newblock {\em Integrability of Poisson
  brackets}. \newblock Math.DG/0210152.

\bibitem{Da}
{\sc Dazord, P.:} \newblock {\em Groupo\"ides 
symplectiques et troisi\`eme th\'eor\`eme de
Lie ``non lin\'eaire''}.
\newblock Lect. Notes in Math. {\bf 1416}, 39--74, Springer-Verlag,
Berlin, 1990.

\bibitem{Duis}
{\sc Duistermaat, J.J.:} \newblock {\em On global action angle coordinates}.
\newblock Comm. Pure Appl. Math. {\bf 33} (1980), 687--706.

\bibitem {GinzLu92} 
{\sc Ginzburg, V.~L., Lu, J.-H.: }\newblock {\em
Poisson cohomology of {M}orita-equivalent {P}oisson manifolds}.
\newblock Internat. Math. Res. Notices {\bf 10} (1992), 199-205.

\bibitem {Gompf95}
{\sc Gompf, R.~E.: }\newblock {\em A new construction of symplectic manifolds}.
\newblock Ann. of Math. (2)  {\bf 142} (1995), 527--595.

\bibitem {Ha84}
{\sc Haefliger, A.:}\newblock {\em  Groupo\"ides d'holonomie et
espaces  classifiants}. 
\newblock Structure transverse des feuilletages. Ast\'erisque {\bf 116}
(1984), 70--97. 

\bibitem {HS87}
{\sc Hilsum, M., Skandalis, G.:}\newblock {\em  Morphismes K-orient\'es d'espaces de feuilles
et fonctorialit\'e en th\'eorie de Kasparov (d'apre\'es une conjecture d'A.Connes)}. 
\newblock Ann. Sci. \'Ecole Norm. Sup. {\bf 20} (1987), 325--390. 

\bibitem {Ho79}
{\sc Howe, R.:}
\newblock {\em $\theta $-series and invariant theory}.
\newblock In: {\em Automorphic forms, representations and $L$-functions
(Proc. Sympos. Pure Math., Oregon State Univ., Corvallis, Ore., 1977),
Part 1}, 275--285.
Proc. Sympos. Pure Math., XXXIII,
Amer. Math. Soc., Providence, R.I., 1979. 


\bibitem {JSW2001}
{\sc Jurco, B., Schupp, P., Wess, J.: }\newblock {\em Noncommutative line
  bundle and Morita equivalence}.
\newblock Lett. Math. Phys. {\bf 61}  (2002), 171--186.

\bibitem{KrMi}
{\sc Kriegl, A.,  Michor, P.W.:} \newblock 
{\em  The Convenient Setting of Global Analysis}.
\newblock {\em Mathematical Surveys and Monographs, 53}.
\newblock American Mathematical Society, Providence, 1997.

\bibitem {Land98}
{\sc Landsman, N.~P.: }\newblock {\em Mathematical Topics between Classical and
  Quantum Mechanics}.
\newblock {\em Springer Monographs in Mathematics}.
\newblock Springer-Verlag, Berlin, Heidelberg, New York, 1998.

\bibitem {Land00b}
{\sc Landsman, N.~P.: }\newblock {\em Bicategories of operator algebras and
  {P}oisson manifolds}.
\newblock In: {\em Mathematical physics in mathematics and physics (Siena,
  2000)}, vol.~30 in {\em Fields Inst. Commun.},   271--286. Amer. Math. Soc.,
  Providence, RI, 2001.


\bibitem{Land01}
{\sc Landsman, N.~P.: } \newblock {\em
 Quantized reduction as a tensor product}
\newblock In: {\em Quantization of
Singular Symplectic Quotients}, N.P. Landsman, M. Pflaum,
M. Schlichenmaier, eds., 137--180, Birkh\"auser, Basel, 2001.

\bibitem {MaXu}
{\sc Mackenzie, K., Xu , P.: }\newblock {\em Integration of {L}ie bialgebroids}.
\newblock Topology  {\bf 39} (2000), 445--467.

\bibitem {Mclan71}
{\sc MacLane, S.: }\newblock {\em Categories for the working mathematician}.
\newblock Springer-Verlag, New York-Berlin, 1971.
\newblock Graduate Texts in Mathematics, Vol. 5.

\bibitem {MiWe}
{\sc Mikami, K., Weinstein, A.: }\newblock {\em Moments and reduction
  for symplectic groupoid actions}.
\newblock Publ. RIMS, Kyoto Univ.  {\bf 24} (1988), 121--140.

\bibitem {Moerd91}
{\sc Moerdijk, I.: }\newblock {\em Classifying toposes and foliations}.
\newblock Ann. Inst. Fourier (Grenoble)  {\bf 41} (1991), 189 -- 289.

\bibitem {Moerd02}
{\sc Moerdijk, I: }\newblock {\em On the classification of regular groupoids}.
\newblock Math.DG/0203099.

\bibitem {MM03}
{\sc Moerdijk, I., Mr\v{c}un., J.: }\newblock {\em Introduction to
  Foliations and Lie Groupoids}, Cambridge University Press (in press).

\bibitem {Morita58}
{\sc Morita, K.: }\newblock {\em Duality for modules and its applications to
  the theory of rings with minimum condition}.
\newblock Sci. Rep. Tokyo Kyoiku Daigaku Sect. A  {\bf 6} (1958), 83--142.

\bibitem {Moser}
{\sc Moser, J.: }\newblock {\em On the volume elements on a manifold}.
\newblock Trans. Amer. Math. Soc.  {\bf 120} (1965), 286--294.

\bibitem {Mrcun99}
{\sc Mr{\v{c}}un, J.: }\newblock {\em Functoriality of the bimodule associated
  to a {H}ilsum-{S}kandalis map}.
\newblock $K$-Theory  {\bf 18} (1999), 235--253.

\bibitem{MuReWi87}
{\sc Muhly, P.S., Renault, J.N., Williams, D.P.:} \newblock {\em Equivalence and
isomorphism for groupoid $C^*$-algebras}.
\newblock J. Operator Theory {\bf 17} (1987), 3--22.

\bibitem{Naga}
{\sc Nagano, T.:} \newblock {\em $1$-forms with the exterior derivative of maximal rank}.
\newblock J. Differential Geom. {\bf 2} (1968), 253--264.

\bibitem{Nielsen}
{\sc Nielsen, J.:} \newblock {\em Untersuchungen zur Theorie der geschlossenen 
zweiseitigen Flachen I}.
\newblock Acta Math.  {\bf 50} (1927), 189--358.

\bibitem{Pr}
{\sc Pronk, D.:} \newblock {\em Etendues and stacks as bicategories of
  fractions}.  \newblock Compositio Math. {\bf 102}  (1996), 243--303.

\bibitem {Rad01}
{\sc Radko, O.: }\newblock {\em A classification of topologically stable
  {P}oisson structures on a compact oriented surface}.
\newblock J. Symp. Geometry. {\bf 1} (2002), 523-542.

\bibitem{Ri74} 
{\sc Rieffel, M.~A.,: }\newblock {\em Morita equivalence for 
$C^{*}$-algebras and $W^{*}$-algebras.} 
\newblock  J. Pure Appl. Algebra {\bf 5}
(1974), 51--96.

\bibitem {RiefSch}
{\sc Rieffel, M.~A., Schwarz, A.: }\newblock {\em Morita equivalence of
  multidimensional noncommutative tori}.
\newblock Internat. J. Math.  {\bf 10} (1999), 289--299.

\bibitem {Severa}
{\sc \v{S}evera, P.:}\newblock {\em Some title 
containing the words ``homotopy'' and ``symplectic'', e.g. this one}.
\newblock Math.SG/0105080.

\bibitem {SeWe01}
{\sc \v{S}evera, P., Weinstein, A.: }\newblock {\em Poisson geometry with a
  $3$-form background}.
\newblock Prog. Theo. Phys. Suppl.  {\bf 144} (2001), 145--154.

\bibitem {TaWe}
{\sc Tang, X.,  Weinstein, A.: } \newblock  {\em Quantization and
  Morita equivalence for 
constant Dirac structures on tori}.  \newblock Math.QA/0305413.

\bibitem{We83} 
{\sc Weinstein, A.: } {\em The local structure of Poisson manifolds}.
\newblock J. Diff. Geom. {\bf 18} (1983), 523--557.

\bibitem {We87}
{\sc Weinstein, A.: }\newblock {\em Symplectic groupoids and {P}oisson
  manifolds}.
\newblock Bull. Amer. Math. Soc. (N.S.)  {\bf 16} (1987), 101--104.

\bibitem {We97}
{\sc Weinstein, A.: }\newblock {\em The modular automorphism group of
  a Poisson manifold}.
\newblock  J. Geom. Phys. {\bf 23} (1997), 379--394.

\bibitem {Wood92}
{\sc Woodhouse, N.~M.~J.: }\newblock {\em Geometric quantization}.
\newblock {\em Oxford Mathematical Monographs}. Second edition.
\newblock The Clarendon Press Oxford University Press, New York, 1992.

\bibitem {Xu91}
{\sc Xu, P.: }\newblock {\em Morita equivalence of {P}oisson manifolds}.
\newblock Comm. Math. Phys.  {\bf 142} (1991), 493--509.

\bibitem {Xu91b}
{\sc Xu, P.: }\newblock {\em Morita equivalent symplectic groupoids}.
\newblock In: {\em Symplectic geometry, groupoids, and integrable systems
  (Berkeley, CA, 1989)},   291--311. Springer, New York, 1991.

\bibitem {Xu}
{\sc Xu, P.: }\newblock {\em Morita equivalence and momentum maps}.
\newblock In preparation .

\bibitem {Zung}
{\sc Zung, N.T.:}\newblock {\em 
Symplectic Topology of Integrable Hamiltonian Systems, II: 
Topological Classification}. \newblock Math.DG/0010181, to appear in
Compositio Math.

\end{thebibliography}
\end{document}